\documentclass[11pt]{article}
\usepackage{amsmath}
\usepackage{amssymb}
\usepackage{amscd}
\usepackage{xspace}
\usepackage{verbatim}
\newcommand{\mysection}[1]{
\section{#1}\setcounter{equation}{0}}
\date{}
\begin{document}
\begin{center}{\bf \Large BOUNDARY VALUE PROBLEMS WITH MEASURES FOR \\[3mm]

\noindent ELLIPTIC EQUATIONS WITH SINGULAR POTENTIALS}
\footnote{Both authors are sponsored by the ECOS-Sud program C08E04. The second author is partially supported by Fondecyt 107125}
 \end{center}
\begin{center}{\bf  Laurent V\'eron}\\
{\small Laboratoire de Math\'ematiques et Physique Th\'eorique}\\
 {\small  Universit\'e Fran\c cois-Rabelais, Tours,  FRANCE} \\[2mm]

 {\bf  Cecilia Yarur}\\
 {\small Departamento de Matem\'aticas y Ciencia de la Computaci\'on} \\
 {\small  Universidad de Santiago de Chile, Santiago, CHILE }
 \end{center}



\newcommand{\txt}[1]{\;\text{ #1 }\;}
\newcommand{\tbf}{\textbf}
\newcommand{\tit}{\textit}
\newcommand{\tsc}{\textsc}
\newcommand{\trm}{\textrm}
\newcommand{\mbf}{\mathbf}
\newcommand{\mrm}{\mathrm}
\newcommand{\bsym}{\boldsymbol}
\newcommand{\scs}{\scriptstyle}
\newcommand{\sss}{\scriptscriptstyle}
\newcommand{\txts}{\textstyle}
\newcommand{\dsps}{\displaystyle}
\newcommand{\fnz}{\footnotesize}
\newcommand{\scz}{\scriptsize}
\newcommand{\be}{
\begin{equation}
}
\newcommand{\bel}[1]{
\begin{equation}
\label{#1}}
\newcommand{\ee}{
\end{equation}
}
\newcommand{\eqnl}[2]{
\begin{equation}
\label{#1}{#2}
\end{equation}
}
\newtheorem{subn}{\name}
\renewcommand{\thesubn}{}
\newcommand{\bsn}[1]{\def\name{#1}
\begin{subn}}
\newcommand{\esn}{
\end{subn}}
\newtheorem{sub}{\name}[section]
\newcommand{\dn}[1]{\def\name{#1}}   
\newcommand{\bs}{
\begin{sub}}
\newcommand{\es}{
\end{sub}}
\newcommand{\bsl}[1]{
\begin{sub}\label{#1}}
\newcommand{\bth}[1]{\def\name{Theorem}
\begin{sub}\label{t:#1}}
\newcommand{\blemma}[1]{\def\name{Lemma}
\begin{sub}\label{l:#1}}
\newcommand{\bcor}[1]{\def\name{Corollary}
\begin{sub}\label{c:#1}}
\newcommand{\bdef}[1]{\def\name{Definition}
\begin{sub}\label{d:#1}}
\newcommand{\bprop}[1]{\def\name{Proposition}
\begin{sub}\label{p:#1}}
\newcommand{\R}{\eqref}
\newcommand{\rth}[1]{Theorem~\ref{t:#1}}
\newcommand{\rlemma}[1]{Lemma~\ref{l:#1}}
\newcommand{\rcor}[1]{Corollary~\ref{c:#1}}
\newcommand{\rdef}[1]{Definition~\ref{d:#1}}
\newcommand{\rprop}[1]{Proposition~\ref{p:#1}} 
\newcommand{\BA}{
\begin{array}}
\newcommand{\EA}{
\end{array}}
\newcommand{\BAN}{\renewcommand{\arraystretch}{1.2}
\setlength{\arraycolsep}{2pt}
\begin{array}}
\newcommand{\BAV}[2]{\renewcommand{\arraystretch}{#1}
\setlength{\arraycolsep}{#2}
\begin{array}}
\newcommand{\BSA}{
\begin{subarray}}
\newcommand{\ESA}{\end{subarray}}
\newcommand{\BAL}{\begin{aligned}}
\newcommand{\EAL}{\end{aligned}}
\newcommand{\BALG}{\begin{alignat}}
\newcommand{\EALG}{\end{alignat}}
\newcommand{\BALGN}{\begin{alignat*}}
\newcommand{\EALGN}{\end{alignat*}}
\newcommand{\note}[1]{\textit{#1.}\hspace{2mm}}
\newcommand{\Proof}{\note{Proof}}
\newcommand{\qeda}{\hspace{10mm}\hfill $\square$}
\newcommand{\qed}{\\
${}$ \hfill $\square$}
\newcommand{\Remark}{\note{Remark}}
\newcommand{\modin}{$\,$\\
[-4mm] \indent}
\newcommand{\forevery}{\quad \forall}
\newcommand{\set}[1]{\{#1\}}
\newcommand{\setdef}[2]{\{\,#1:\,#2\,\}}
\newcommand{\setm}[2]{\{\,#1\mid #2\,\}}
\newcommand{\lra}{\longrightarrow}
\newcommand{\sgn}{\rm{sgn}}
\newcommand{\lla}{\longleftarrow}
\newcommand{\llra}{\longleftrightarrow}
\newcommand{\Lra}{\Longrightarrow}
\newcommand{\Lla}{\Longleftarrow}
\newcommand{\Llra}{\Longleftrightarrow}
\newcommand{\warrow}{\rightharpoonup}
\newcommand{
\paran}[1]{\left (#1 \right )}
\newcommand{\sqbr}[1]{\left [#1 \right ]}
\newcommand{\curlybr}[1]{\left \{#1 \right \}}
\newcommand{\abs}[1]{\left |#1\right |}
\newcommand{\norm}[1]{\left \|#1\right \|}
\newcommand{
\paranb}[1]{\big (#1 \big )}
\newcommand{\lsqbrb}[1]{\big [#1 \big ]}
\newcommand{\lcurlybrb}[1]{\big \{#1 \big \}}
\newcommand{\absb}[1]{\big |#1\big |}
\newcommand{\normb}[1]{\big \|#1\big \|}
\newcommand{
\paranB}[1]{\Big (#1 \Big )}
\newcommand{\absB}[1]{\Big |#1\Big |}
\newcommand{\normB}[1]{\Big \|#1\Big \|}

\newcommand{\thkl}{\rule[-.5mm]{.3mm}{3mm}}
\newcommand{\thknorm}[1]{\thkl #1 \thkl\,}
\newcommand{\trinorm}[1]{|\!|\!| #1 |\!|\!|\,}
\newcommand{\bang}[1]{\langle #1 \rangle}
\def\angb<#1>{\langle #1 \rangle}
\newcommand{\vstrut}[1]{\rule{0mm}{#1}}
\newcommand{\rec}[1]{\frac{1}{#1}}
\newcommand{\opname}[1]{\mbox{\rm #1}\,}
\newcommand{\supp}{\opname{supp}}
\newcommand{\dist}{\opname{dist}}
\newcommand{\myfrac}[2]{{\displaystyle \frac{#1}{#2} }}
\newcommand{\myint}[2]{{\displaystyle \int_{#1}^{#2}}}
\newcommand{\mysum}[2]{{\displaystyle \sum_{#1}^{#2}}}
\newcommand {\dint}{{\displaystyle \int\!\!\int}}
\newcommand{\q}{\quad}
\newcommand{\qq}{\qquad}
\newcommand{\hsp}[1]{\hspace{#1mm}}
\newcommand{\vsp}[1]{\vspace{#1mm}}
\newcommand{\ity}{\infty}
\newcommand{\prt}{
\partial}
\newcommand{\sms}{\setminus}
\newcommand{\ems}{\emptyset}
\newcommand{\ti}{\times}
\newcommand{\pr}{^\prime}
\newcommand{\ppr}{^{\prime\prime}}
\newcommand{\tl}{\tilde}
\newcommand{\sbs}{\subset}
\newcommand{\sbeq}{\subseteq}
\newcommand{\nind}{\noindent}
\newcommand{\ind}{\indent}
\newcommand{\ovl}{\overline}
\newcommand{\unl}{\underline}
\newcommand{\nin}{\not\in}
\newcommand{\pfrac}[2]{\genfrac{(}{)}{}{}{#1}{#2}}

\def\ga{\alpha}     \def\gb{\beta}       \def\gg{\gamma}
\def\gc{\chi}       \def\gd{\delta}      \def\ge{\epsilon}
\def\gth{\theta}                         \def\vge{\varepsilon}
\def\gf{\phi}       \def\vgf{\varphi}    \def\gh{\eta}
\def\gi{\iota}      \def\gk{\kappa}      \def\gl{\lambda}
\def\gm{\mu}        \def\gn{\nu}         \def\gp{\pi}
\def\vgp{\varpi}    \def\gr{\rho}        \def\vgr{\varrho}
\def\gs{\sigma}     \def\vgs{\varsigma}  \def\gt{\tau}
\def\gu{\upsilon}   \def\gv{\vartheta}   \def\gw{\omega}
\def\gx{\xi}        \def\gy{\psi}        \def\gz{\zeta}
\def\Gg{\Gamma}     \def\Gd{\Delta}      \def\Gf{\Phi}
\def\Gth{\Theta}
\def\Gl{\Lambda}    \def\Gs{\Sigma}      \def\Gp{\Pi}
\def\Gw{\Omega}     \def\Gx{\Xi}         \def\Gy{\Psi}

\def\CS{{\mathcal S}}   \def\CM{{\mathcal M}}   \def\CN{{\mathcal N}}
\def\CR{{\mathcal R}}   \def\CO{{\mathcal O}}   \def\CP{{\mathcal P}}
\def\CA{{\mathcal A}}   \def\CB{{\mathcal B}}   \def\CC{{\mathcal C}}
\def\CD{{\mathcal D}}   \def\CE{{\mathcal E}}   \def\CF{{\mathcal F}}
\def\CG{{\mathcal G}}   \def\CH{{\mathcal H}}   \def\CI{{\mathcal I}}
\def\CJ{{\mathcal J}}   \def\CK{{\mathcal K}}   \def\CL{{\mathcal L}}
\def\CT{{\mathcal T}}   \def\CU{{\mathcal U}}   \def\CV{{\mathcal V}}
\def\CZ{{\mathcal Z}}   \def\CX{{\mathcal X}}   \def\CY{{\mathcal Y}}
\def\CW{{\mathcal W}} \def\CQ{{\mathcal Q}} 
\def\BBA {\mathbb A}   \def\BBb {\mathbb B}    \def\BBC {\mathbb C}
\def\BBD {\mathbb D}   \def\BBE {\mathbb E}    \def\BBF {\mathbb F}
\def\BBG {\mathbb G}   \def\BBH {\mathbb H}    \def\BBI {\mathbb I}
\def\BBJ {\mathbb J}   \def\BBK {\mathbb K}    \def\BBL {\mathbb L}
\def\BBM {\mathbb M}   \def\BBN {\mathbb N}    \def\BBO {\mathbb O}
\def\BBP {\mathbb P}   \def\BBR {\mathbb R}    \def\BBS {\mathbb S}
\def\BBT {\mathbb T}   \def\BBU {\mathbb U}    \def\BBV {\mathbb V}
\def\BBW {\mathbb W}   \def\BBX {\mathbb X}    \def\BBY {\mathbb Y}
\def\BBZ {\mathbb Z}

\def\GTA {\mathfrak A}   \def\GTB {\mathfrak B}    \def\GTC {\mathfrak C}
\def\GTD {\mathfrak D}   \def\GTE {\mathfrak E}    \def\GTF {\mathfrak F}
\def\GTG {\mathfrak G}   \def\GTH {\mathfrak H}    \def\GTI {\mathfrak I}
\def\GTJ {\mathfrak J}   \def\GTK {\mathfrak K}    \def\GTL {\mathfrak L}
\def\GTM {\mathfrak M}   \def\GTN {\mathfrak N}    \def\GTO {\mathfrak O}
\def\GTP {\mathfrak P}   \def\GTR {\mathfrak R}    \def\GTS {\mathfrak S}
\def\GTT {\mathfrak T}   \def\GTU {\mathfrak U}    \def\GTV {\mathfrak V}
\def\GTW {\mathfrak W}   \def\GTX {\mathfrak X}    \def\GTY {\mathfrak Y}
\def\GTZ {\mathfrak Z}   \def\GTQ {\mathfrak Q}

\font\Sym= msam10 
\def\SYM#1{\hbox{\Sym #1}}
\newcommand{\bdw}{\prt\Gw\xspace}
\medskip
\begin{abstract}
We study the boundary value problem with Radon measures  for nonnegative
solutions of $-\Delta u+Vu=0$ in a bounded smooth domain
$\Gw$, when $V$ is a locally bounded nonnegative function. Introducing some specific capacity, we give sufficient conditions on a Radon measure $\gm$ on $\prt\Gw$ so that the problem can be solved. We study the reduced measure associated to this equation as well as the boundary trace of positive solutions.  
\end{abstract}

\noindent
{\it \footnotesize 1991 Mathematics Subject Classification}. {\scriptsize
35J25; 28A12; 31C15; 47B44; 35C15}.\\
{\it \footnotesize Key words}. {\scriptsize Laplacian; Poisson
potential; capacities; singularities; Borel measures; Harnack inequalities. }
\tableofcontents
\mysection{Introduction}
 Let $\Gw$ be a smooth bounded domain of $\BBR^N$ and $V$ a locally bounded real valued measurable function defined in $\Gw$. The first question we adress is the solvability of the following non-homogeneous Dirichlet problem with a Radon measure for boundary data,
\begin{equation}\label{I1}
\left\{\BA{ll}
-\Gd u+Vu=0\qquad&\text{in }\Gw\\
\phantom{-\Gd u+V}
u=\gm\qquad&\text{in }\prt\Gw.
\EA\right.
\end{equation}
Let $\gr$ be the first (and positive) eigenfunction of $-\Gd$ in $W^{1,2}_0(\Gw)$. By a solution we mean a function $u\in L^1(\Gw)$, such that $Vu\in L^1_\gr$, which satisfies
 \begin{equation}\label{I2}
\myint{\Gw}{}\left(-u\Gd\gz+Vu\gz\right)dx=-\myint{\prt\Gw}{}\myfrac{\prt\gz}{\prt{\bf n}}d\gm.
\end{equation}
for any function $\gz\in C_0^1(\overline\Gw)$ such that $\Gd\gz\in L^\infty(\Gw)$. When $V$ is a bounded nonnegative function, it is straightforward that there exist a unique solution. However, it is less obvious to find general conditions which allow the solvability for any $\gm\in \frak M(\prt\Gw)$, the set of Radon measures on $\prt\Gw$. In order to avoid difficulties due to Fredholm type obstructions, we shall most often assume that $V$ is nonnegative, in which case there exists at most one solution. \smallskip

Let us denote by $K^\Gw$ the Poisson kernel in $\Gw$ and by $\BBK[\gm]$ the Poisson potential of a measure, that is
 \begin{equation}\label{I6}
\BBK[\gm](x):=\myint{\prt\Gw}{}K^\Gw(x,y)d\gm(y)\qquad\forall x\in\Gw.
\end{equation}
We first observe that, when $V\geq 0$ and the measure $\gm$ satisfies
 \begin{equation}\label{I5}
\myint{\Gw}{}\BBK[|\gm|](x)V(x)\gr(x)dx<\infty, 
\end{equation}
then problem $(\ref{I1})$ admits a solution. A Radon measure which satisfies $(\ref{I5})$ is called {\it an admissible measure} and a measure for which a solution exists is called {\it a good measure}. \smallskip

We first consider the {\it subcritical case} which means that the boundary value is solvable for any  $\gm\in \frak M(\prt\Gw)$. As a first result, we prove that any measure $\gm$ is admissible if $V$ is nonnegative and satisfies
 \begin{equation}\label{I4}
\sup_{\;\;\;\;\;\;y\in\prt\Gw}\!\!\!\!\! {\rm ess}\myint{\Gw}{}K^\Gw(x,y)V(x)\gr(x)dx<\infty.
\end{equation}
Using estimates on the Poisson kernel, this condition is fulfilled if there exists $M>0$ such that for any $y\in\prt\Gw$,
 \begin{equation}\label{I3}
\myint{0}{D(\Gw)}\left(\myint{\Gw\cap B_r(y)}{}V(x)\gr^2(x) dx\right)\myfrac{dr}{r^{N+1}}\leq M
\end{equation}
where $D(\Gw)=diam (\Gw)$. We give also sufficient conditions which ensures that the boundary value problem $(\ref{I1})$ is stable from the weak*-topology of $\frak M(\prt\Gw)$ to $L^1(\Gw)\cap L^1_{V\gr}(\Gw)$. One of the sufficient conditions is that  $V\geq 0$ satisfies
\begin{equation}\label{stab0}
\lim_{\ge\to 0}\myint{0}{\ge}\left(\myint{\Gw\cap B_r(y)}{}V(x)\gr^2(x)dx\right)\myfrac{dr}{r^{N+1}}=0,
\end{equation} 
uniformly with respect to $y\in\prt\Gw$.
\smallskip

In the {\it supercritical case} problem $(\ref{I1})$ cannot be solved for any $\gm\in \frak M(\prt\Gw)$. In order to characterize positive good measures, we introduce a framework of nonlinear analysis which have been used by Dynkin and Kuznetsov (see \cite{Dbook1} and references therein) and Marcus and V\'eron \cite{MV1}  in their study of the  boundary value problems with measures
\begin{equation}\label{I8}
\left\{\BA{ll}
-\Gd u+|u|^{q-1}u=0\qquad&\text{in }\Gw\\
\phantom{-\Gd u+|u|^{q-1}}
u=\gm\qquad&\text{in }\prt\Gw,
\EA\right.
\end{equation}
where $q>1$. In these works, positive good measures on $\prt\Gw$ are completely characterized by the $C_{2/q,q'}$-Bessel in dimension N-1 and the following property:

\smallskip {\it A measure $\gm\in\frak M_+(\prt\Gw)$ is good for problem $(\ref{I8})$ if and only if it does charge Borel sets with zero $C_{2/q,q'}$-capacity, i.e
\begin{equation}\label{I9}
C_{2/q,q'}(E)=0\Longrightarrow\gm(E)=0\qquad\forall E\subset\prt\Gw,\,E\text{ Borel}.
\end{equation}
Moreover, any positive good measure is the limit of an increasing sequence $\{\gm_n\}$ of admissible measures which, in this case, are the positive measures belonging to the Besov space $B_{2/q,q'}(\prt\Gw)$.
They also characaterize  removable sets in terms of $C_{2/q,q'}$-capacity.}\medskip

In our present work, and always with $V\geq 0$, we use a capacity associated to the Poisson kernel $K^\Gw$ and belongs to a class  studied by Fuglede \cite{Fu} \cite{Fu2}. It is defined by
\begin{equation}\label{I10}
C_V(E)=\sup\{\gm (E):\gm\in\frak M_+(\prt\Gw), \gm(E^c)=0,\,\norm{V\BBK[\gm]}_{L^1_\gr}\leq 1\},
\end{equation}
for any  Borel set $E\subset\prt\Gw$. Furtheremore $C_V(E)$ is equal to the value of its dual expression $C^*_V(E)$ defined by
\begin{equation}\label{I11}
C^*_V(E)=\inf\{\norm {f}_{L^\infty}:\check\BBK[f]\geq 1\quad\text {on }E\},
\end{equation}
where
\begin{equation}\label{I12}
\check\BBK[f](y)=\myint{\Gw}{}K^\Gw(x,y)f(x)V(x)\gr(x)dx\qquad\forall y\in\prt\Gw.
\end{equation}
If $E$ is a compact subset of $\prt\Gw$, this capacity is explicitely given by
\begin{equation}\label{I12'}
C_V(E)=C^*_V(E)=\max_{y\in E}\left(\myint{\Gw}{}K^\Gw(x,y)V(x)\gr(x) dx\right)^{-1}.
\end{equation}
\medskip

We denote by $Z_V$ the largest set with zero $C_V$ capacity, i.e.

\begin{equation}\label{I7}
Z_V=\left\{y\in \prt\Gw:\myint{\Gw}{}K^\Gw(x,y)V(x)\gr(x) dx=\infty\right\},
\end{equation}
and we prove the following.\smallskip 

\noindent{\it 1- If $\{\gm_n\}$ is an increasing sequence of positive good measures which converges to a measure $\gm$  in the weak* topology, then $\gm$ is a good measure. 
\smallskip 

\noindent 2- If $\gm\in \frak M_+(\prt\Gw)$ satisfies $\gm (Z_V)=0$, then $\gm$ is a good measure.\smallskip 

\noindent 3- A good measure $\gm$ vanishes on $Z_V$ if and only if there exists  an increasing sequence of positive admissible measures which converges to $\gm$  in the weak* topology.
}\medskip

In section 4 we study relaxation phenomenon in replacing $(\ref{I1})$ by the truncated problem
\begin{equation}\label{I1-k}
\left\{\BA{ll}
-\Gd u+V_ku=0\qquad&\text{in }\Gw\\
\phantom{-\Gd u+V_k}
u=\gm\qquad&\text{in }\prt\Gw.
\EA\right.
\end{equation}
where $\{V_{k}\}$ is an increasing sequence of positive bounded functions which converges to $V$
locally uniformly in $\Gw$. We adapt to the linear problem some of the principles of the reduced measure. This notion is introduced by Brezis, Marcus and Ponce \cite{BMP} in the study of the nonlinear Poisson equation
\begin{equation}\label{I10'}
-\Gd u+g(u)=\gm\qquad\text{in }\Gw\
\end{equation}
and extended to the Dirichlet problem 
\begin{equation}\label{I10''}
\left\{\BA{ll}
-\Gd u+g(u)=0\qquad&\text{in }\Gw\\
\phantom{-\Gd u+g()}
u=\gm\qquad&\text{in }\prt\Gw,
\EA\right.
\end{equation}
by Brezis and Ponce \cite{BP}.  In our construction, problem $(\ref{I1-k})$ admits a unique solution $u_k$. The sequence $\{u_k\}$ decreases and converges to some $u$ which satisfies a relaxed boundary value problem
\begin{equation}\label{I13}
\left\{\BA{ll}
-\Gd u+Vu=0\qquad&\text{in }\Gw\\
\phantom{-\Gd u+V}
u=\gm^*\qquad&\text{in }\prt\Gw.
\EA\right.
\end{equation}
The measure $\gm^*$ is called the {\it reduced measure} associated to $\gm$ and $V$. {\it Note that $\gm^*$ is the largest measure for which the problem
\begin{equation}\label{I15}
\left\{\BA{ll}
-\Gd u+Vu=0\qquad&\text{in }\Gw\\
\phantom{-\Gd u+V}
u=\gn\leq \gm\qquad&\text{in }\prt\Gw.
\EA\right.
\end{equation}
admits a solution}. This truncation process allows to construct the Poisson kernel $K^\Gw_V$ associated to the operator $-\Gd +V$ as being the limit of the decreasing limit of the sequence of kernel functions $\{K^\Gw_{V_k}\}$ asociated to $-\Gd +V_k$.
The solution $u=u_{\gm^*}$ of $(\ref{I13})$ is expressed by
\begin{equation}\label{I14}
u_{\gm^*}(x)=\myint{\prt\Gw}{}K^\Gw_V(x,y)d\gm(y)=\myint{\prt\Gw}{}K^\Gw_V(x,y)d\gm^*(y)\qquad\forall x\in\Gw.
\end{equation}
We define the vanishing set of $K_V$ by
\begin{equation}\label{I16}
Z^*_V=\{y\in \prt\Gw:K^\Gw_V(x_0,y)=0\},
\end{equation}
for some $x_0\in \Gw$, and thus for any $x\in \Gw$ by Harnack inequality. We prove 
\smallskip

\noindent 1- $Z^*_V\subset Z_V$.
\smallskip

\noindent 2- $\gm^*=\gm\chi_{_{Z^*_V}}$.\smallskip

\noindent A challenging open problem is to give conditions on $V$ which allows $Z^*_V= Z_V$. \smallskip

The last section is devoted to the construction of the boundary trace of positive solutions of 
\begin{equation}\label{I19}
-\Gd u+Vu=0\qquad\text{in }\Gw,
\end{equation}
assuming $V\geq 0$. Using results of \cite{MV3}, we defined the regular set $\CR(u)$ of the boundary trace of $u$. This set is a relatively open subset of $\prt\Gw$ and the regular part of the boundary trace is represented by a positive Radon measure $\gm_u$ on $\CR(u)$. In order to study the singular set of the boundary trace $\CS(u):=\prt\Gw\setminus \CR(u)$, we adapt the sweeping method introduced by Marcus and V\'eron in \cite{MV4} for equation 
\begin{equation}\label{I20}
-\Gd u+g(u)=0\qquad\text{in }\Gw.
\end{equation}
If $\gm$ is a good positive measure concentrated on $\CS(u)$, and $u_\gm$ is the unique solution of $(\ref{I1})$ with boundary data $\gm$, we set $v_\gm=\min\{u,u_\gm\}$. Then $v_\gm$ is a positive super solution which admits a positive trace $\gg_u(\gm)\in\frak M_+(\prt\Gw)$. The extended boundary trace $Tr^e(u)$ of $u$ is defined by
\begin{equation}\label{I21}
\gn(u)(E):=Tr^e(u)(E)=\sup\{\gg_u(\gm)(E):\gm\text { good}, E\subset\prt\Gw,\,E\text{ Borel}\}.
\end{equation}
 {\it Then $Tr^e(u)$ is a Borel measure on $\Gw$. If we assume moreover that 
 \begin{equation}\label{I22}
\lim_{\ge\to 0}\myint{0}{\ge}\left(\myint{\Gw\cap B_r(y)}{}V(x)\gr^2(x)dx\right)\myfrac{dr}{r^{N+1}}=0\qquad\text{uniformly with respect to }y\in\prt\Gw,
\end{equation}
then $Tr^e(u)$ is a bounded measure and therefore a Radon measure. Finally, if $N=2$ and $(\ref{I22})$ holds, or if $N=2$ and there holds
 \begin{equation}\label{I22'}
\lim_{\ge\to 0}\myint{0}{\ge}\left(\myint{\Gw\cap B_r(y)}{}V(x)(\gr(x)-\ge)^2_+dx\right)\myfrac{dr}{r^{N+1}}=0,
\end{equation}
uniformly with respect to $\ge\in (0,\ge_0] \text{ and }y\text{ s.t. }\dist(x,\prt\Gw)=\ge$, then $u=u_{\gn(u)}$.}\medskip

If $V(x)\leq v(\gr(x)$ for some $v$ which satisfies 
 \begin{equation}\label{I23}
\myint{0}{1}v(t)tdt<\infty,
\end{equation}
then Marcus and V\'eron proved in \cite{MV3} that $u=u_{\gn_u}$. Actually, when $V$ has such a geometric form, the assumptions $(\ref{I22})$-$(\ref{I22'})$ and $(\ref{I23})$ are equivalent.
\mysection{The subcritical case}
In the sequel $\Gw$ is a bounded smooth domain in $\BBR^N$ and $V\in L^\infty_{loc}$. We denote by $\gr$ the first eigenfunction of $-\Gd$ in $W^{1,2}_0(\Gw)$, $\gr>0$ with the corresponding eigenvalue $\gl$, by
$\frak M(\prt\Gw)$ the space of bounded Radon measures on $\prt\Gw$ and by $\frak M_+(\prt\Gw)$ its positive cone. For any positive Radon measure on $\prt\Gw$, we shall denote by the same symbol the corresponding outer regular bounded Borel measure. Conversely, for any outer regular bounded Borel $\gm$, we denote by the same expression $\gm$ the Radon measure defined on $C(\prt\Gw)$ by
$$\gz\mapsto \gm(\gz)=\myint{\prt\Gw}{}\gz d\gm.
$$
If $\gm\in \frak M(\prt\Gw)$, we are concerned with the following problem
\begin{equation}\label{bvp1}
\left\{\BA{ll}
-\Gd u+Vu=0\qquad&\text{in }\Gw\\
\phantom{-\Gd u+V}
u=\gm\qquad&\text{in }\prt\Gw.
\EA\right.
\end{equation}

\bdef{bvpdef} Let  $\gm\in \frak M(\prt\Gw)$. We say that $u$ is a weak solution of $(\ref{bvp1})$, if $u\in L^1(\Gw)$, $Vu\in L^1_\gr(\Gw)$ and, for any $\gz\in C^1_0(\overline\Gw)$ with $\Gd \gz\in L^\infty(\Gw)$, there holds
 \begin{equation}\label{bvp2}
\myint{\Gw}{}\left(-u\Gd\gz+Vu\gz\right)dx=-\myint{\prt\Gw}{}\myfrac{\prt\gz}{\prt {\bf n}}d\gm.
\end{equation}
\es
In the sequel we put
$$T(\Gw):=\{\gz\in C^1_0(\overline\Gw)\text { such that }\Gd \gz\in L^\infty(\Gw)\}.
$$

We recall the following estimates obtained by Brezis \cite{Br2}
\bprop{brez} Let  $\gm\in L^1(\prt\Gw)$ and $u$ be a weak solution of problem $(\ref{bvp1})$. Then there holds
 \begin{equation}\label{brez1}
\norm u_{L^1(\Gw)}+\norm {V_+u}_{L_\gr^1(\Gw)}\leq \norm {V_{-}u}_{L_\gr^1(\Gw)}+c\norm \gm_{L^1(\prt\Gw)}
\end{equation}
 \begin{equation}\label{brez2}
\myint{\Gw}{}\left(-|u|\Gd\gz+V|u|\gz\right)dx\leq-\myint{\prt\Gw}{}\myfrac{\prt\gz}{\prt {\bf n}}|\gm|dS
\end{equation}
and 
 \begin{equation}\label{brez3}
\myint{\Gw}{}\left(-u_+\Gd\gz+Vu_+\gz\right)dx\leq-\myint{\prt\Gw}{}\myfrac{\prt\gz}{\prt {\bf n}}\gm_+dS,
\end{equation}
for all $\gz\in T(\Gw)$, $\gz\geq 0$.
\es

We denote by $K^{\Gw}(x,y)$ the Poisson kernel in $\Gw$ and by $\BBK[\gm]$ the Poisson potential of $\gm\in \frak M(\prt\Gw)$ defined by
 \begin{equation}\label{Pois}
\BBK[\gm](x)=\myint{\prt\Gw}{}K^{\Gw}(x,y)d\gm(y)\qquad\forall x\in\Gw.
\end{equation}

\bdef {adm} A measure $\gm$ on $\prt\Gw$ is {\bf admissible} if
 \begin{equation}\label{adm}
\myint{\Gw}{}\BBK[|\gm|](x)|V(x)|\gr(x)dx<\infty.
\end{equation}
It is {\bf good} if problem $(\ref{bvp1})$ admits a weak solution.
\es

We notice that, if there exists at least one admissible positive measure $\gm$, then
\begin{equation}\label{glob}
\myint{\Gw}{}V(x)\gr^2(x)dx<\infty.
\end{equation}
\bth {Exist} Assume $V\geq 0$, then problem $(\ref{bvp1})$ admits at most one solution.  Furthermore, if $\gm$ is admissible, then there exists a unique solution that we denote $u_\gm$.
\es
\Proof Uniqueness follows from $(\ref{brez1})$. For existence we can assume $\gm\geq 0$. For any $k\in\BBN_{*}$ set $V_{k}=\inf\{V,k\}$ and denote by $u:=u_{k}$ the solution of 
\begin{equation}\label{L4}\left\{\BA {ll}
-\Gd u+V_{k}(x)u=0\qquad&\text {in }\Gw\\
\phantom{-\Gd u+V(x)}
u=\gm\qquad&\text {on }\prt\Gw.
\EA\right.\end{equation} 
Then $0\leq u_{k}\leq \BBK[\gm]$. 
By the maximum principle, $u_{k}$ is decreasing and converges to some $u$, and 
$$0\leq V_{k}u_{k}\leq V\BBK[\gm]. 
$$
Thus, by dominated convergence theorem $V_{k}u_{k}\to Vu$ in $L^1_\gr$. Setting $\gz\in T(\Gw)$ and letting $k$ tend to infinity in  equality
\begin{equation}\label{L5}
\myint{\Gw}{}\left(-u_{k}\Gd\gz+V_{k}u_{k}\gz\right)dx=-\myint{\prt\Gw}{}\frac{\prt\gz}{\prt{\bf n}}d\gm, 
\end{equation}
implies that $u$ satisfies $(\ref{bvp2})$.\qeda\medskip

\noindent\Remark If $V$ changes sign, we can put $\tilde u=u+\BBK[\gm]$. Then $(\ref{bvp1})$ is equivalent to
\begin{equation}\label{bvp4}
\left\{\BA{ll}
-\Gd \tilde u+V\tilde u=V\BBK[\gm]&\qquad\text{in }\Gw\\
\phantom{-\Gd \tilde u+V}
\tilde u=0\qquad&\text{in }\prt\Gw.
\EA\right.
\end{equation}
This is a Fredholm type problem (at least if the operator $\phi\mapsto R(v):=(-\Gd)^{-1}(V\phi)$ is compact in $L^1_\gr(\Gw)$). Existence will be ensured by orthogonality conditions. \medskip

If we assume that $V\geq 0$ and
\begin{equation}\label{bvp5}
\myint{\Gw}{}K^\Gw(x,y)V(x)\gr(x)dx<\infty,
\end{equation}
for some $y\in \prt\Gw$, then $\gd_y$ is admissible. The following result yields to the  solvability of $(\ref{bvp1})$ for any $\gm\in\frak M_+(\Gw)$.
\bprop{Uncond} Assume $V\geq 0$ and the integrals $(\ref{bvp5})$ are bounded uniformly with respect to $y\in\prt\Gw$. Then any measure on $\prt\Gw$ is admissible.
\es
\Proof If $M$ is the upper bound of these integrals and $\gm\in \frak M_+(\prt\Gw)$, we have,
\begin{equation}\label{ener4}
\myint{\Gw}{}\BBK[\gm](x)V(x)\gr(x)dx=\myint{\prt\Gw}{}\left(\myint{\Gw}{}K^\Gw(x,y)V(x)\gr(x)dx\right)d\gm(y)\leq 
M\gm(\prt\Gw),
\end{equation}
 by Fubini's theorem. Thus $\gm$ is admissible.
\qeda\medskip

\noindent\Remark Since the Poisson kernel in $\Gw$ satisfies the two-sided estimate
\begin{equation}\label{ener7}
c^{-1}\myfrac{\gr(x)}{|x-y|^{N}}\leq K^\Gw(x,y)\leq c\myfrac{\gr(x)}{|x-y|^{N}}\qquad \forall (x,y)\in \Gw\ti\prt\Gw,
\end{equation}
for some $c>0$, assumption $(\ref{bvp5})$ is equivalent to 
\begin{equation}\label{ener6}
\myint{\Gw}{}\myfrac{V(x)\gr^2(x)}{|x-y|^{N}}dx<\infty.
\end{equation}
This implies $(\ref{glob})$ in particular. If we set $D_y=\max\{|x-y|:x\in\Gw\}$, then

$$\BA {l}\myint{\Gw}{}\myfrac{V(x)\gr^2(x)}{|x-y|^{N}}dx= 
\myint{0}{D_y}\left(\myint{\{x\in\Gw:|x-y|=r\}}{}\!\!\!\!\!\!V(x)\gr^2(x)dS_r(x)\right)\myfrac{dr}{r^{N}}\\[4mm]
\phantom{\myint{\Gw}{}\myfrac{V(x)\gr(x)}{|x-y|^{N-1}}dx}
=\displaystyle\lim_{\ge\to 0}\left(\left[r^{-N}\myint{\Gw\cap B_r(y)}{}\!\!\!\!\!\!V(x)\gr^2(x)dx\right]_\ge^{D_y}
+N\myint{\ge}{D_y}\left(\myint{\Gw\cap B_r(y)}{}\!\!\!\!\!\!V(x)\gr^2(x)dx\right)\myfrac{dr}{r^{N+1}}\right)
\EA$$
(both quantity may be infinite). Thus, if we assume 
 \begin{equation}\label{marc0}
\myint{0}{D_y}\left(\myint{\Gw\cap B_r(y)}{}\!\!\!\!\!\!V(x)\gr^2(x)dx\right)\myfrac{dr}{r^{N+1}}<\infty,
 \end{equation}
there holds 
 \begin{equation}\label{Marc0'}
 \liminf_{\ge\to 0}\ge^{-N}\myint{\Gw\cap B_\ge(y)}{}\!\!\!\!\!\!V(x)\gr^2(x)dS=0.
 \end{equation}
Consequently
 \begin{equation}\label{Marc0''}
\myint{\Gw}{}\myfrac{V(x)\gr^2(x)}{|x-y|^{N}}dx=
D_y^{-N}\myint{\Gw}{}V(x)\gr^2(x)dx+N\myint{0}{D_y}\left(\myint{\Gw\cap B_r(y)}{}\!\!\!\!\!\!V(x)\gr^2(x)dx\right)\myfrac{dr}{r^{N+1}}.
 \end{equation}
Therefore $(\ref{bvp5})$ holds and $\gd_y$ is admissible. \medskip


As a natural extension of \rprop{Uncond}, we have the following stability result.
\bth{stab} Assume $V\geq 0$ and 
\begin{equation}\label{cond'}
\lim_{\tiny\BA {l}E\text{ Borel}\\|E|\to 0\EA}\myint{E}{}K^\Gw(x,y)V(x)\gr(x)dx=0\quad\text {uniformly with respect to }y\in\prt\Gw.
\end{equation}
If $\gm_n$ is a sequence of positive Radon measures on $\prt\Gw$ converging to $\gm$ in the weak* topology, then $u_{\gm_n}$ converges to $u_{\gm}$ in $L^1(\Gw)\cap L^1_{V\gr} (\Gw)$ and locally uniformly in $\Gw$.
\es
\Proof 
We put $u_{\gm_n}:=u_n$. By the maximum principle $0\leq u_{n}\leq \BBK[\gm_n]$. Furthermore, it follows from $(\ref{brez1})$ that
 \begin{equation}\label{stab1}
\norm {u_{n}}_{L^1(\Gw)}+\norm {Vu_{n}}_{L_\gr^1(\Gw)}\leq c\norm {\gm_n}_{L^1(\prt\Gw)}\leq C.
\end{equation}
Since $-\Gd u_n$ is bounded in $L_\gr^1(\Gw)$, the sequence $\{u_n\}$ is relatively compact in $L^1(\Gw)$ by the regularity theory for elliptic equations. Therefore, there exist a subsequence  $u_{n_k}$ and some function $u\in L^1(\Gw)$ with $Vu\in L^1_\gr(\Gw)$ such that
$u_{n_k}$ converges to $u$ in $L^1(\Gw)$, almost everywhere on $\Gw$  and locally uniformly in $\Gw$ since $V\in L^\infty_{loc}(\Gw)$. The main question is to prove the convergence of $Vu_{n_k}$ in $L^1_\gr(\Gw)$. If $E\subset\Gw$ is any Borel set, there holds
$$\BA {l}
\myint{E}{}u_{n}V(x)\gr(x) dx\leq\myint{E}{}\BBK[\gm_n]V(x)\gr(x) dx\\[3mm]\phantom{\myint{E}{}u_{n}V(x)\gr(x) dx}
\leq\myint{\prt\Gw}{}\left(\myint{E}{}K^\Gw(x,y)V(x)\gr(x) dx\right)d\gm_n(y)\\[3mm]\phantom{\myint{E}{}u_{n}V(x)\gr(x) dx}
\leq M_n\displaystyle\max_{y\in\prt\Gw}\myint{E}{}K^\Gw(x,y)V(x)\gr(x) dx
,\EA$$
where $M_n:=\gm_n(\prt\Gw)$. Thus
 \begin{equation}\label{stab5}\BA {l}
\myint{E}{}u_{n}V(x)\gr(x) dx
\leq M_n\displaystyle\max_{y\in\prt\Gw}\myint{E}{}K^\Gw(x,y)V(x)\gr(x) dx.
\EA\end{equation}
Then, by $(\ref{cond'})$,
$$
\lim_{|E|\to 0}\myint{E}{}u_{n}V(x)\gr(x) dx=0.
$$
As a consequence the set of function $\{u_{n}\gr V\}$ is uniformly integrable. By Vitali's theorem
$Vu_{n_k}\to Vu$ in $L^1_\gr(\Gw)$. Since 
 \begin{equation}\label{stab6}
\myint{\Gw}{}\left(-u_n\Gd\gz+Vu_n\gz\right)dx=-\myint{\prt\Gw}{}\myfrac{\prt\gz}{\prt {\bf n}}d\gm_n,
\end{equation}
for any $\gz\in T(\Gw)$, the function $u$ satisfies $(\ref{bvp2})$.
\qeda\medskip

Assumption $(\ref{cond'})$ may be difficult to verify and the following result gives an easier formulation.

\bprop{CS} Assume $V\geq 0$ satisfies
 \begin{equation}\label{Marc0}
\lim_{\ge\to 0}\myint{0}{\ge}\left(\myint{\Gw\cap B_r(y)}{}V(x)\gr^2(x)dx\right)\myfrac{dr}{r^{N+1}}=0\quad\text{uniformly with respect to }y\in\prt\Gw.
\end{equation} 
Then $(\ref{cond'})$ holds.
\es
\Proof If $E\subset\Gw$ is a Borel set and $\gd>0$, we put $E_\gd=E\cap B_\gd(y)$ and $E^c_\gd=E\setminus E_\gd$. Then
$$ \myint{E}{}\myfrac{V(x)\gr^2(x)}{|x-y|^N}dx= \myint{E_\gd}{}\myfrac{V(x)\gr^2(x)}{|x-y|^N}dx+ \myint{E_\gd^c}{}\myfrac{V(x)\gr^2(x)}{|x-y|^N}dx.
$$
Clearly
 \begin{equation}\label{Marc1}
 \myint{E_\gd^c}{}\myfrac{V(x)\gr^2(x)}{|x-y|^N}dx\leq \gd^{-N} \myint{E.}{}V(x)\gr^2(x) dx.
\end{equation} 
Since  $(\ref{marc0})$ holds for any $y\in\prt\Gw$, $(\ref{Marc0''})$ implies
 \begin{equation}\label{Marc3}
\myint{E_\gd}{}\myfrac{V(x)\gr^2(x)}{|x-y|^N}dx=
\gd^{-N}\myint{E_\gd}{}V(x)\gr^2(x)dx+N\myint{0}{\gd}\left(\myint{E\cap B_r(y)}{}V(x)\gr^2(x)dx\right)\myfrac{dr}{r^{N+1}}.
\end{equation} 
Using $(\ref{Marc0})$, for any $\ge>0$, there exists $s_0>0$ such that for any  $s>0$ and $y\in\prt\Gw$
$$s\leq s_0\Longrightarrow N\myint{0}{s}\left(\myint{B_r(y)}{}V(x)\gr^2(x)dx\right)\myfrac{dr}{r^{N+1}}\leq\ge/2.
$$
We fix $\gd=s_0$. Since $(\ref{glob})$ holds,
\begin{equation}\label{cond+}
\lim_{\tiny\BA {l}E\text{ Borel}\\|E|\to 0\EA}\myint{E}{}V(x)\gr^2(x)dx=0.
\end{equation}
Then there exists $\eta>0$ such that for any Borel set $E\subset\Gw$,
$$|E|\leq\eta\Longrightarrow \myint{E}{}V(x)\gr^2(x)dx\leq s^N_0\ge/4.
$$
Thus
$$\myint{E}{}\myfrac{V(x)\gr^2(x)}{|x-y|^N}dx\leq\ge.
$$
This implies the claim by $(\ref{ener7})$.\qeda\medskip

An assumption which is used in \cite[Lemma 7.4]{MV3} in order to prove the existence of a boundary trace of any positive solution of $(\ref{I19})$ is that there exists some nonnegative measurable function $v$ defined on $\BBR_+$ such that
 \begin{equation}\label{t1}
\abs {V(x)}\leq v(\gr(x))\quad\forall x\in\Gw\quad\text{and }\myint{0}{s}tv(t)dt<\infty\quad\forall s>0.
\end{equation}

In the next result we show that condition $(\ref{t1})$ implies $(\ref{cond'})$.

\bprop{Tr} Assume $V$ satisfies $(\ref{t1})$. Then
\begin{equation}\label{condabs}
\lim_{\tiny\BA {l}E\text{ Borel}\\|E|\to 0\EA}\myint{E}{}K^\Gw(x,y)\abs{V(x)}\gr(x)dx=0\quad\text {uniformly with respect to }y\in\prt\Gw.
\end{equation}
\es
\Proof Since $\prt\Gw$ is $C^2$, there exist $\ge_0>0$ such that any for any $x\in\Gw$ satisfying $\gr(x)\leq\ge_0$, there exists a unique $\gs(x)\in \prt\Gw$ such that $|x-\gs(x)|=\gr(x)$. 
We use $(\ref{Marc0})$ in \rprop{CS} under the equivalent form
 \begin{equation}\label{Marc0+}
\lim_{\ge\to 0}\myint{0}{\ge}\left(\myint{\Gw\cap C_r(y)}{}|V(x)|\gr^2(x)dx\right)\myfrac{dr}{r^{N+1}}=0\quad\text{uniformly with respect to }y\in\prt\Gw,
\end{equation} 
in which we have replaced $ B_r(y)$ by the the cylinder $C_r(y):=\{x\in\Gw:\gr(x)<r,|\gs(x)-y|<r\}$. Then
$$\BA {l}
\myint{0}{\ge}\left(\myint{\Gw\cap C_r(y)}{}|V(x)|\gr^2(x)dx\right)\myfrac{dr}{r^{N+1}}
\leq c\myint{0}{\ge}\left(\myint{0}{r}v(t)t^2dt\right)\myfrac{dr}{r^{2}}\\[4mm]
\phantom{\myint{0}{\ge}\left(\myint{\Gw\cap C_r(y)}{}|V(x)|\gr^2(x)dx\right)\myfrac{dr}{r^{N+1}}}
\leq  c\myint{0}{\ge}v(t)\left(1-\myfrac{t}{\ge}\right)tdt\\[4mm]
\phantom{\myint{0}{\ge}\left(\myint{\Gw\cap C_r(y)}{}|V(x)|\gr^2(x)dx\right)\myfrac{dr}{r^{N+1}}}
\leq c\myint{0}{\ge}v(t)tdt.
\EA$$
Thus $(\ref{Marc0})$ holds.\qeda\medskip


\mysection{The capacitary approach}

Throughout this section $V$ is a locally bounded nonnegative and measurable function defined on $\Gw$. We assume that there exists a positive measure $\gm_0$ on $\prt\Gw$ such that 
\begin{equation}\label{capa1}
\myint{\Gw}{}\BBK[\gm_0]V(x)\gr(x)dx=\CE(1,\gm_0)<\infty.
\end{equation}

\bdef {energy}  If $\gm\in\frak M_+(\prt\Gw)$ and $f$ is a nonnegative measurable function defined in $\Gw$ such that 
$$(x,y)\mapsto \BBK[\gm](y)f(x)V(x)\gr(x)\in L^1(\Gw\ti\prt\Gw;dx\otimes d\gm),$$
we set
\begin{equation}\label{capa2}
\CE(f,\gm)=\myint{\Gw}{}\left(\myint{\prt\Gw}{}K^\Gw(x,y)d\gm(y)\right)f(x)V(x)\gr(x)dx.
\end{equation}
\es
If we put
\begin{equation}\label{capa3}
\check\BBK_V[f](y)=\myint{\Gw}{}K^\Gw(x,y)f(x)V(x)\gr(x)dx,
\end{equation}
then, by Fubini's theorem, $\check\BBK_{V}[f]<\infty$, $\gm$-almost everywhere on $\prt\Gw$ and
\begin{equation}\label{capa4}
\CE(f,\gm)=\myint{\prt\Gw}{}\left(\myint{\Gw}{}K^\Gw(x,y)f(x)V(x)\gr(x)dx\right)d\gm(y).
\end{equation}
\bprop {lsc} Let $f$ be fixed. Then \smallskip

\noindent (a) $y\mapsto \check \BBK_V[f](y)$ is lower semicontinuous on $\prt\Gw$.\smallskip

\noindent (b) $\gm\mapsto \CE(f,\gm)$ is lower semicontinuous on $\frak M_+(\prt\Gw)$ in the weak*-topology
\es
\Proof Since $y\mapsto K^\Gw(x,y)$ is continuous, statement (a) follows by Fatou's lemma. If ${\gm_n}$ is a sequence in $\frak M_+(\prt\Gw)$ converging to some $\gm$ in the weak*-topology, then 
$\BBK[\gm_n]$ converges to $\BBK[\gm]$ everywhere in $\Gw$. By Fatou's lemma 
$$\CE(f,\gm)\leq \liminf_{n\to\infty}\myint{\Gw}{}\BBK[\gm_n](x)f(x)V(x)\gr(x)dx=\liminf_{n\to\infty}\CE(f,\gm_n).
$$
\qeda

Notice that if $V\gr f\in L^p(\Gw)$, for $p>N$, then $\BBG[Vf\gr]\in C^1(\overline\Gw)$ and
 \begin{equation}\label{capa5}
\check \BBK[f](y):=\myint{\Gw}{}K^\Gw(x,y)V(x)f(x)\gr(x)dx=-\myfrac{\prt}{\prt{\bf n}}\BBG[Vf\gr](y).
\end{equation}
This is in particular the case if $f$ has compact support in $\Gw$.\medskip

\bdef{mspace} We denote by $\frak M^V(\prt\Gw)$ the set of all measures $\gm$ on $\prt\Gw$ such that
$V\BBK[\gm]\in L^1_{\gr}(\Gw)$. If $\gm$ is such a measure, we denote
 \begin{equation}\label{mspace1}
\norm\gm_{\frak M^V}= \myint{\Gw}{}\abs{\BBK[\gm](x)}V(x)\gr(x) dx=\norm{V\BBK[\gm]}_{L^1_\gr}.
\end{equation}
\es

Clearly $\norm{\,.\,}_{\frak M^V}$ is a norm. The space  $\frak M^V(\prt\Gw)$ is not complete but its positive cone $\frak M_+^V(\prt\Gw)$ is complete. If $E\subset\prt\Gw$ is a Borel subset, we put
$$\frak M_+(E)=\{\gm\in\frak M_+(\prt\Gw): \gm (E^c)=0\}\quad\text{and }\;\frak M^V_+(E)=\frak M_+(E)\cap \frak M^V(\prt\Gw).
$$
\bdef{cap} If $E\subset\prt\Gw$ is any Borel subset we set
 \begin{equation}\label{capa6}
C_V(E):=\sup \{\gm (E) :\gm\in\frak M^V_+(E), \norm\gm_{\frak M^V}\leq 1\}.
\end{equation}
\es
We notice that $(\ref{capa6})$ is equivalent to
 \begin{equation}\label{fug1}
C_V(E):=\sup \left\{\myfrac{\gm (E)}{\norm\gm_{\frak M^V}} :\gm\in\frak M^V_+(E)\right\}.
\end{equation}

\bprop{capaci} The set function $C_V$ satisfies. 
 \begin{equation}\label{fu0}
C_V(E)\leq\sup_{y\in E} \left(\myint{\Gw}{}K^\Gw(x,y)V(x)\gr(x)dx\right)^{-1}\quad\forall E\subset\prt\Gw,\, E\text{ Borel},
\end{equation}
and equality holds in $(\ref{fu0})$ if $E$ is compact. Moreover, 
 \begin{equation}\label{fug0}
C_V(E_1\cup E_2)=\sup \{C_V(E_1),C_V(E_2)\}\quad\forall E_i\subset\prt\Gw,\, E_i\text{ Borel}.
\end{equation}
\es
\Proof Notice that $E\mapsto C_V(E)$ is a nondecreasing set function for the inclusion relation and that $(\ref{capa6})$ implies
 \begin{equation}\label{fu1}
\gm (E)\leq C_V(E)\norm\gm_{\frak M^V}\qquad\forall\gm\in\frak M^V_+(E).
\end{equation}
Let $E\subset\prt\Gw$ be a Borel set and $\gm\in \frak M_+(E)$. Then
$$\BA {l}\norm {\gm}_{\frak M^V}=\myint{E}{}\left(\myint{\Gw}{}K^{\Gw}(x,y)V(x)\gr(x) dx\right)d\gm(y)\\[4mm]
\phantom{\norm {\gm_{\frak M^V}}}\geq \gm(E)\displaystyle{\inf_{y\in E}\myint{\Gw}{}K^{\Gw}(x,y)V(x)\gr(x) dx}.
\EA$$
Using $(\ref{capa6})$ we derive
 \begin{equation}\label{fu2}
C_V(E)\leq \sup_{y\in E}\left(\myint{\Gw}{}K^{\Gw}(x,y)V(x)\gr(x) dx\right)^{-1}.
\end{equation}
If $E$ is compact, there exists $y_0\in E$ such that
$$\inf_{y\in E}\myint{\Gw}{}K^{\Gw}(x,y)V(x)\gr(x) dx=\myint{\Gw}{}K^{\Gw}(x,y_0)V(x)\gr(x) dx,
$$
since $y\mapsto \check \BBK[1](y)$ is l.s.c.. Thus
$$\norm {\gd_{y_0}}_{\frak M^V}=\gd_{y_0}(E)\myint{\Gw}{}K^{\Gw}(x,y_0)V(x)\gr(x) dx
$$
and 
$$C_V(E)\geq \myfrac{\gd_{y_0}(E)}{\norm {\gd_{y_0}}_{\frak M^V}}=\sup_{y\in E}\left(\myint{\Gw}{}K^{\Gw}(x,y)V(x)\gr(x) dx\right)^{-1}.
$$
Therefore equality holds in $(\ref{fu0})$.
Identity $(\ref{fug0})$ follows  $(\ref{fu0})$ when there is equality. Moreover it holds if $E_1$ and $E_2$ are two arbitrary compact sets. Since $C_V$ is eventually an inner regular capacity (i.e. $C_V(E)=\sup \{C_V(K): K\subset E,\, K\text { compact}\}$) it holds for any Borel set. However we give below a self-contained proof.
If  $E_1$ and $E_2$ be two disjoint Borel subsets of $\prt\Gw$, for any $\ge>0$ 
there exists $\gm\in \frak M^V_+(E_1\cup E_2)$ such that 
$$\myfrac{\gm (E_1)+\gm (E_2)}{\norm{\gm}_{\frak M^V}}\leq C_V(E_{1}\cup E_{2})\leq \myfrac{\gm (E_1)+\gm (E_2)}{\norm{\gm}_{\frak M^V}}+\ge.
$$
Set $\gm_i=\chi_{_{E_i}}\gm$. Then $\gm_i\in \frak M^V_+(E_i)$ and  $\norm{\gm}_{\frak M^V}=\norm{\gm_1}_{\frak M^V}+\norm{\gm_2}_{\frak M^V}$. By $(\ref{fu1})$ 
 \begin{equation}\label{capa6*}
 C_V(E_{1}\cup E_{2})\leq \myfrac{\norm{\gm_1}_{\frak M^V}}{\norm{\gm_1}_{\frak M^V}+\norm{\gm_2}_{\frak M^V}}C_V(E_{1})+\myfrac{\norm{\gm_2}_{\frak M^V}}{\norm{\gm_1}_{\frak M^V}+\norm{\gm_2}_{\frak M^V}}C_V(E_{2})+\ge
\end{equation}
This implies that there exists $\gth\in [0,1]$ such that
 \begin{equation}\label{capa6**}
 C_V(E_{1}\cup E_{2})\leq \gth C_V(E_{1})+(1-\gth)C_V(E_{2})\leq \max\{C_V(E_{1}),C_V(E_{2})\}.
\end{equation}
Since $C_V(E_{1}\cup E_{2})\geq\max\{C_V(E_{1}),C_V(E_{2})\}$ as $C_V$ is increasing,
 \begin{equation}\label{capa6***}
  E_1\cap E_2=\emptyset\Longrightarrow C_V(E_{1}\cup E_{2})=\max\{C_V(E_{1}),C_V(E_{2})\}.
\end{equation}
If  $E_1\cap E_2\neq\emptyset$, then $E_1\cup E_2=E_1\cup (E_2\cap E^c_1)$
and therefore
$$C_V(E_{1}\cup E_{2})=\max\{C_V(E_{1}),C_V(E_2\cap E^c_1)\}\leq \max\{C_V(E_{1}),C_V(E_{2})\}.
$$
Using again $(\ref{fug1})$ we derive $(\ref{fug0})$.
\qeda\medskip

The following set function is the dual expression of $C_V(E)$.

\bdef {dual}For any Borel set $E\subset\prt\Gw$, we set
 \begin{equation}\label{capa7}
C^*_V(E):=\inf \{\norm f_{L^\infty} :\check \BBK[f](y)\geq 1\quad\forall y\in E\}.
\end{equation}
\es 
The next result is stated in \cite[p 922]{Fu2} using minimax theorem and the fact that $K^\Gw$ is lower semi continuous in $\Gw\ti\prt\Gw$. Although the proof is not explicited, a simple adaptation of the proof of \cite[Th 2.5.1]{AH} leads to the result. 

\bprop{Equi} For any compact set  $E\subset\prt\Gw$,
\begin{equation}\label {zeroset}
C_V(E)= C^*_V(E).
 \end{equation}
 \es

In the same paper \cite{Fu2}, formula $(\ref{fu0})$ with equality is claimed (if $E$ is compact).

\bth{G} If $\{\gm_{n}\}$ is an increasing sequence of good measures converging to some measure $\gm$ in the weak* topology, then $\gm$ is good.
\es
\Proof We use formulation $(\ref{D9})$. We take for test function the function $\eta$ solution of
\begin{equation}\left\{\BA {ll}
-\Gd \eta=1\qquad&\text{in }\Gw\\
\phantom{-\Gd}\eta=0\qquad&\text{on }\Gw,
\EA\right.\end{equation}
there holds
$$\myint{\Gw}{}\left(1+V\right)u_{\gm_{n}}\eta dx=
-\myint{\prt\Gw}{}\myfrac{\prt\eta}{\prt{\bf n}}d\gm_{n}\leq c^{-1}\gm_{n}(\prt\Gw)\leq c^{-1}\gm(\prt\Gw)
$$
where $c>0$ is such that
$$c^{-1}\geq -\myfrac{\prt\eta}{\prt{\bf n}}\geq c\quad\text {on }\prt\Gw.
$$
Since $\{u_{\gm_{n}}\}$ is increasing and $\eta\leq c\gr$ by Hopf boundary lemma, we can let $n\to\infty$ by the monotone convergence theorem. If $u:=\lim_{n\to\infty}u_{\gm_{n}}$, we obtain
$$\myint{\Gw}{}\left(1+V\right)u \eta dx\leq c^{-1}\gm(\prt\Gw).
$$
Thus $u$ and $\gr Vu$ are in $L^1(\Gw)$. Next, if 
$\gz\in C_{0}^{1}(\overline\Gw)\cap C^{1,1}(\overline\Gw)$, then
$u_{\gm_{n}}|\Gd\gz|\leq Cu_{\gm_{n}}$ and $Vu_{\gm_{n}}|\gz|\leq CVu_{\gm_{n}}\eta$. Because the sequence $\{u_{\gm_{n}}\}$ and 
$\{Vu_{\gm_{n}}\eta\}$ are uniformly integrable, the same holds for 
$\{u_{\gm_{n}}\Gd\gz\}$ and $\{Vu_{\gm_{n}}\gz\}$. Considering
$$
\myint{\Gw}{}\left(-u_{\gm_{n}}\Gd\gz+Vu_{\gm_{n}}\gz\right)dx=-\myint{\prt\Gw}{}\myfrac{\prt\gz}{\prt{\bf n}}d\gm_{n}.
$$
it follows by Vitali's theorem, 
$$
\myint{\Gw}{}\left(-u\Gd\gz+Vu\gz\right)dx=-\myint{\prt\Gw}{}\myfrac{\prt\gz}{\prt{\bf n}}d\gm.
$$
Thus $\gm$ is a good measure.\qeda\medskip

We define the {\it singular boundary set} $Z_V$ by
\begin{equation}\label{Z0}
Z_V=\left\{y\in\prt\Gw:\myint{\Gw}{}K^\Gw(x,y)V(x)\gr(x) dx=\infty\right\}.
\end{equation}
Since $\check\BBK[1]$ is l.s.c., it is a Borel function and $Z_V$ is a Borel set. The next result characterizes  the good measures.
\bprop{van} Let $\gm$ be an admissible positive measure. Then $\gm(Z_V)=0$.
\es
\Proof If $K\subset Z_V$ is compact,  $\gm_K=\chi_{_K}\gm$ is admissible, thus, by Fubini theorem
$$\norm{\gm_K}_{\frak M^V}=\myint{K}{}\left(\myint{\Gw}{}K^{\Gw}(x,y)V(x)\gr(x)dx\right)d\gm(y)<\infty.
$$
Since 
$$\myint{\Gw}{}K^{\Gw}(x,y)V(x)\gr(x)dx\equiv \infty\qquad\forall y\in K
$$
it follows that $\gm(K)=0$. This implies $\gm(Z_V)=0$ by regularity.\qeda
\bth{approx} Let $\gm\in \frak M_+(\prt\Gw)$ such that
\begin{equation}\label{Z1}
\gm(Z_V)=0.
\end{equation}
Then $\gm $ is good.
\es
\Proof Since $\check\BBK[1]$ is l.s.c., for any $n\in\BBN_*$, 
$$K_n:=\{y\in \prt\Gw:\check\BBK[1](y)\leq n\}
$$
is a compact subset of $\prt\Gw$. Furthermore $K_n\cap Z_V=\emptyset$ and $\cup K_n=Z_V^c$. Let 
$\gm_n=\chi_{_{K_n}}\gm$, then
\begin{equation}\label{Z4}
\CE(1,\gm_n)=\myint{\Gw}{}\BBK[\gm_n]V(x)\gr(x)dx\leq n\gm_n(K_n).
\end{equation}
Therefore $\gm_n$ is admissible. By the monotone convergence theorem, $\gm_n\uparrow \chi_{_{Z_{V^c}}}\gm$ and by \rth{G}, $\chi_{_{Z_{V^c}}}\gm$ is good. Since $(\ref{Z1})$ holds, $\chi_{_{Z_{V^c}}}\gm=\gm$, which ends the proof.
\qeda\medskip

The full characterization of the good measures in the general case appears to be difficult without any further assumptions on $V$. However the following  holds
\bth{equiv} Let $\gm\in \frak M_+(\prt\Gw)$ be a good measure. The following assertions are equivalent:\smallskip

\noindent (i) $\gm(Z_V)=0$.\smallskip

\noindent (ii) There exists an increasing sequence of admissible measures $\{\gm_n\}$ which converges to $\gm$ in the weak*-topology.
\es
\Proof If (i) holds, it follows from the proof of \rth{approx} that the sequence $\{\gm_n\}$ increases and converges to $\gm$. If (ii) holds, any admissible measure $\gm_n$ vanishes on $Z_V$ by \rprop{van}. Since
$\gm_n\leq \gm$, there exists an increasing sequence of $\gm$-integrable functions $h_n$ such that $\gm_n=h_n \gm$. Then
$\gm_n(Z_V)$ increases to $\gm(Z_V)$ by the monotone convergence theorem. The conclusion follows from the fact that $\gm_n(Z_V)=0$. \qeda


\mysection{Representation formula and reduced measures}

We recall the construction of the Poisson kernel for $-\Gd +V$: if we look for a solution of 
 \begin{equation}\label{D1}\left\{
 \BA {ll}
 -\Gd v+V(x)v=0\qquad&\text {in }\Gw\\ 
 \phantom{ -\Gd v+V(x)}
 v=\gn\qquad&\text {in }\prt\Gw,
 \EA\right.
\end{equation}
where $\gn\in\frak M(\prt\Gw)$, $V\geq 0$, $V\in L^\infty_{loc}(\Gw)$, we can consider an increasing sequence of smooth domains 
$\Gw_n$ such that $\overline\Gw_n\subset\Gw_{n+1}$ and $\cup_n\Gw_n=\cup_n\overline\Gw_n=\Gw$. For each of these domains, denote by $K^{\Gw}_{V\chi_{_{\Gw_n}}}$ the Poisson kernel of $-\Gd+V\chi_{_{\Gw_n}}$ in $\Gw$ and by $\BBK_{V\chi_{_{\Gw_n}}}[.]$ the corresponding operator. We denote by $K^\Gw:=K^{\Gw}_0$ the Poisson kernel in $\Gw$ and by $\BBK[.]$ the Poisson operator in $\Gw$. Then the solution $v:=v_n$ of
 \begin{equation}\label{D2}\left\{
 \BA {ll}
 -\Gd v+V\chi_{_{\Gw_n}}v=0\qquad&\text {in }\Gw\\ 
 \phantom{  -\Gd v+V\chi_{_{\Gw_n}}}
 v=\gn&\text {in }\prt\Gw,
 \EA\right.
\end{equation}
is expressed by
 \begin{equation}\label{D3}
 v_n(x)=\myint{\prt\Gw}{}K^{\Gw}_{V\chi_{_{\Gw_n}}}(x,y)d\gn (y)=\BBK_{V\chi_{_{\Gw_n}}}[\gn](x).
\end{equation}
If $G^\Gw$ is the Green kernel of $-\Gd$ in $\Gw$ and $\BBG[.]$ the corresponding Green operator, $(\ref{D3})$ is equivalent to 
 \begin{equation}\label{D4}
v_n(x)+\myint{\Gw}{}G^\Gw(x,y)(V\chi_{_{\Gw_n}}v_n)(y)dy=\myint{\prt\Gw}{}K^{\Gw}(x,y)d\gn (y),
\end{equation}
equivalently
$$v_n+\BBG[V\chi_{_{\Gw_n}}v_n]=\BBK[\gn].
$$
Notice that this equality is equivalent to the weak formulation of problem $(\ref{D2})$: for any $\gz\in T(\Gw)$, there holds
 \begin{equation}\label{D4'}
\myint{\Gw}{}\left(-v_n\Gd\gz+V\chi_{_{\Gw_n}}v_n\gz\right)dx=-\myint{\prt\Gw}{}\myfrac{\prt\gz}{\prt{\bf n}}d\gn.
\end{equation}
Since $n\mapsto K^{\Gw}_{V\chi_{_{\Gw_n}}}$ is decreasing, the sequence $\{v_n\}$ inherits this property and there exists 
 \begin{equation}\label{D5}
\lim_{n\to\infty}K^{\Gw}_{V\chi_{_{\Gw_n}}}(x,y)=K^{\Gw}_{V}(x,y).
\end{equation}
By the monotone convergence theorem, 
 \begin{equation}\label{D6}  
\lim_{n\to\infty}v_n(x)=v(x)=\myint{\prt\Gw}{}K^{\Gw}_{V}(x,y)d\gn(y).
 \end{equation}
By Fatou's theorem
 \begin{equation}\label{D7}
 \myint{\Gw}{}G^\Gw(x,y)V(y)v(y)dy\leq \liminf_{n\to\infty}\myint{\Gw}{}G^\Gw(x,y)(V\chi_{_{\Gw_n}}v_n)(y)dy,
\end{equation}
and thus,  
 \begin{equation}\label{D8}
v(x)+ \myint{\Gw}{}G^\Gw(x,y)V(y)v(y)dy\leq\BBK[\gn](x)\qquad\forall x\in \Gw.
\end{equation}
Now the main question is to know whether $v$ keeps the boundary value $\gn$. Equivalently, whether the equality holds in $(\ref{D7})$ with $\lim$ instead of $\liminf$, and therefore in $(\ref{D8})$. This question is associated to the notion of reduced measured in the sense of Brezis-Marcus-Ponce: Since $Vv\in L^1_{\gr}(\Gw)$ and 
 \begin{equation}\label{D9}
 -\Gd v+V(x)v=0\qquad\text {in }\Gw
 \end{equation}
 holds, the function $v+\BBG[Vv]$ is positive and harmonic in $\Gw$. Thus it admits a boundary trace $\gn^*\in \frak M_{+}(\prt\Gw)$ and
 \begin{equation}\label{D10}
v+ \BBG[Vv]=\BBK[\gn^*].
\end{equation}
Equivalently $v$  satisfies the relaxed problem
 \begin{equation}\label{D11}\left\{
 \BA {ll}
 -\Gd v+V(x)v=0&\qquad\text {in }\Gw\\ 
 \phantom{ -\Gd v+V(x)}
 v=\gn^*\qquad&\text {in }\prt\Gw,
 \EA\right.
 \end{equation}
 and thus $v=u_{\gn^*}$. Noticed that $\gn^*\leq\gn$ and the mapping $\gn\mapsto\gn^*$ is nondecreasing.
\bdef{red}The measure $\gn^*$ is the {\it reduced measure} associated to $\gn$.\es
\bprop {opt} There holds $\BBK_V[\gn]=\BBK_V[\gn^*]$. Furthermore
the reduced measure $\gn^*$ is the largest measure  for which the following problem  
\begin{equation}\label{D12}\left\{
 \BA {ll}
  \phantom{.} -\Gd v+V(x)v=0\qquad&\text {in }\Gw\\ 
  \phantom{}
\!\gl\in \frak M_{+}(\prt\Gw),\;\gl\leq\gn
 \\ 
 \phantom{ -\Gd v+V(x);.}
 v=\gl\qquad&\text {in }\prt\Gw,
 \EA\right.
 \end{equation}
 admits a solution.
\es
\Proof The first assertion follows from the fact that $v=\BBK_V[\gn]$ by $(\ref{D5})$ and $v=u_{\gn^*}=\BBK_V[\gn^*]$ by $(\ref{D11})$.
It is clear that $\gn^*\leq\gn$ and that the problem $(\ref{D12})$ admits a solution for $\gl=\gn^*$. If $\gl$ is a positive measure smaller than 
$\gm$, then $\gl^*\leq \gm^*$. But if there exist some $\gl$ such that the problem $(\ref{D12})$ admits a solution, then $\gl=\gl^*$. This implies the claim.\qeda\medskip

As a consequence of the characterization of $\gn^*$ there holds

\bcor{Uni} Assume $V\geq 0$ and let $\{V_k\}$ be an increasing sequence of nonnegative bounded measurable functions converging to $V$ a.e. in $\Gw$. Then the solution $u_k$ of 
 \begin{equation}\label{uni1}\left\{\BA {ll}
-\Gd u+V_ku=0\qquad&\text{in }\Gw\\ \phantom{-\Gd u+V_k}
u=\gn\qquad&\text{in }\prt\Gw,
\EA\right.\end{equation}
converges to $u_{\gn^*}$.
\es
\Proof The previous construction shows that 
$u_k=\BBK_{V_k}[\gn]$ decreases to some $\tilde u$ which satisfies a relaxed equation, the boundary data of which, $\tilde\gn^*$,  is the largest measure $\gl\leq\gn$ for which problem $(\ref{D12})$ admits a solution. Therefore $\tilde\gn^*=\gn^*$ and $\tilde u=u_{\gn^*}$. Similarly  $\{K^\Gw_{V_k}\}$  decreases and converges to $K^\Gw_{V}$. \qeda\medskip

We define the boundary vanishing set of $K^\Gw_V$  by

 \begin{equation}\label{van1}
Z^*_V:=\{y\in\prt\Gw\,|\,K^\Gw_V(x,y)=0\}\quad\text{for some } x\in\Gw.
\end{equation}
Since $V\in L^\infty_{loc}(\Gw)$, $Z^*_V$ is independent of $x$ by Harnack inequality; furthermore it is a Borel set.

\bth{vanish}  Let $\gn\in\frak M_+(\prt\Gw)$.\smallskip

\noindent (i) If $\gn((Z^*_V)^c)=0$, then $\gn^*=0$.\smallskip

\noindent (ii) There always holds $Z^*_V\subset Z_V$. 
\es
\Proof The first assertion is clear since $\gn=\chi_{_{Z^*_V}}\gn+\chi_{_{(Z^*_V)^c}}\gn=\chi_{_{Z^*_V}}\gn$ and, by \rprop{opt},
$$u_{\gn^*}(x)=\BBK_V[\gn^*](x)=\myint{Z^*_V}{}K^\Gw_V(x,y)d\gn(y)=0\qquad\forall x\in\Gw,
$$
by definition of $Z^*_V$. For proving (ii), we assume that $C_V( Z^*_V)>0$; there exists $\gm\in\frak M^V_+( Z^*_V)$ such that $\gm(Z^*_V)>0$. 
Since $\gm$ is admissible let $u_\gm$ be the solution of $(\ref{I1})$. Then $\gm^*=\gm$, thus $u_\gm=\BBK^V[\gm]$
and
$$\BBK^V[\gm](x)=\myint{\prt\Gw}{}K^\Gw_V(x,y)d\gm(y)=\myint{Z^*_V}{}K^\Gw_V(x,y)d\gm(y)=0,
$$
contradiction. Thus $C_V(Z^*_V)=0$. Since $(\ref{fu0})$  implies that $Z_V$ is the largest Borel set with zero $C_V$-capacity, it implies $Z^*_V\subset Z _V$.\qeda



\medskip

In order to obtain more precise informations on $Z^*_V$ some minimal regularity assumptions on $V$ are needed. We also recall the following result proved by Ancona \cite{An}.
\bth{An} Assume $V\geq 0$ satisfies $\gr^2V\in L^\infty(\Gw)$. If for some $y_0\in \prt\Gw$ and any cone 
$C_{y_0}$ with vertex $y_0$ having the property that 
$\overline {C}_{y_0}\cap B_r(y_0)\subset \Gw\cup\{y_0\}$ for some $r>0$, there exists $c_1>0$ such that 
 \begin{equation}\label{D14}
\forall (x,y)\in \Gw\cap B_r(y_0)\ti  \Gw\cap B_r(y_0), \,|x-y_0|=|y-y_0|\leq r\Longrightarrow 
c^{-1}\leq\myfrac{V(x)}{V(y)}\leq c_1
\end{equation}
and
 \begin{equation}\label{D15}
 \myint{0}{r}V(t{\bf n}_{y_0})tdt=\infty,
\end{equation}
 where $\bf {n}_0$ is the normal outward unit vector to $\prt\Gw$ at $y_0$, then
 \begin{equation}\label{D16}
K_V^{\Gw}(x,y_0)=0\qquad\forall x\in \Gw.
\end{equation}
\es

We define the {\it conical  singular boundary set}
 \begin{equation}\label{D17}
\tilde Z_V=\left\{y\in \prt\Gw:\myint{\Gw\cap C_y}{}K^\Gw(x,y)V(x)\gr(x) dx=\infty\;\text{for some cone }C_y\Subset\Gw\right\}
\end{equation}
where $C_y\Subset\Gw$ means that there exists $a>0$ such that $\overline C_y\cap B_a(y)\subset \Gw\cup\{y\}$. Clearly $\tilde Z_V\subset  Z_V$. 

\bcor{Osc} Assume $V\geq 0$ satisfies $\gr^2V\in L^\infty(\Gw)$ and the conical oscillation condition $(\ref{D14})$ of \rth{An} for any $y\in Z_V$. Then $\tilde Z_V=Z^*_V$.
\es
\Proof We can assume that $y=0$ and denote $C_y=C$. Since 
$$K^\Gw(x,0)V(x)\gr(x)\leq ca^{-N}V(x)\gr^2(x)\qquad\forall x\in \Gw\cap B_a^c,$$
and $V\gr^2\in L^1(\Gw)$, there holds, using $(\ref{ener7})$,
$$\myint{B_a\cap C}{}V(x)\gr^2(x) \myfrac{dx}{|x|^N}=\infty.
$$
Using spherical coordinates and the fact that $\gr^2(x)\geq c|x|$ in $B_a\cap C_y$,
$$\myint{0}{a}\myint{S}{}V(r,\gs)r d\gs\,dr=\infty.
$$
where $S=C\cap \prt B_1$. But in $C\cap B_a$ the oscillation condition $(\ref{D14})$ holds. This implies
 \begin{equation}\label{D18}
 \myint{0}{a}V(r,\gs)tdt=\infty\qquad\forall \gs\in S.
\end{equation}
Thus $y\in Z^*_V$.\qeda


\mysection{The boundary trace}

\subsection{The regular part}
In this section, $V\in L^\infty_{loc}(\Gw)$ is nonnegative. If $0<\ge\leq \ge_0$, we denote $d(x)=\dist (x,\prt\Gw)$ for $x\in\Gw$, and set $\Gw_\ge:=\{x\in\Gw:d(x)>\ge\}$, $\Gw'_\ge=\Gw\setminus\Gw_\ge$ and $\Gs_\ge=\prt\Gw_\ge$. It is well known that there exists $\ge_0$ such that, for any $0<\ge\leq \ge_0$  and any $x\in\Gw'_\ge$ there exists a unique projection $\gs(x)$ of $x$ on $\prt\Gw$ and any $x\in \Gw'_\ge$ can be written in a unique way under the form
$$x=\gs(x)-d(x){\bf n}
$$
where $\bf n$ is the outward normal unit vector to $\prt\Gw$ at $\gs(x)$.
The mapping $x\mapsto(d(x),\gs(x))$ is a $C^2$ diffeomorphism from $\Gw'_\ge$ to $(0,\ge_0]\ti\prt\Gw$. We recall the following definition given in \cite{MV3}. If $\CA$ is a Borel subset of $\prt\Gw$, we set 
$\CA_\ge=\{x\in \Gs_\ge:\gs(x)\in A\}$.
\bdef{Trdef} Let $\CA$ be a relatively open subset of $\prt\Gw$, $\{\gm_\ge\}$ be a set of Radon measures on $\CA_{\ge}$ $(0<\ge\leq \ge_0)$ and $\gm\in\frak M(\CA)$. We say that $\gm_\ge\rightharpoonup\gm$ in the weak*-topology if, for any $\gz\in C_c(\CA)$, 
 \begin{equation}\label{Y1}
\lim_{\ge\to 0}\myint{\CA_{\ge}}{}\gz(\gs(x))d\gm_\ge(x)=\myint{\CA}{}\gz d\gm.
\end{equation}
A function $u\in C(\Gw)$ possesses a boundary trace $\gm\in \frak M(\CA)$ if
 \begin{equation}\label{Y2}
\lim_{\ge\to 0}\myint{\CA_{\ge}}{}\gz(\gs(x))u(x)dS(x)=\myint{\CA}{}\gz d\gm\qq\forall \gz\in C_c(\CA).
\end{equation}
\es
The following result is proved in \cite[p 694]{MV3}.
\bprop{Trlem1} Let $u\in C(\Gw)$ be a positive solution of 
 \begin{equation}\label{F1}
-\Gd u+V(x)u=0\qquad\text{in }\Gw.
\end{equation}
Assume that, for some $z\in\prt\Gw$, there exists an open neighborhood $U$ of $z$ such that 
 \begin{equation}\label{F2}
\myint{U\cap\Gw}{}Vu\gr(x)dx<\infty.
\end{equation}
Then $u\in L^1(K\cap\Gw)$ for any compact subset $K\subset G$ and there exists a positive Radon measure $\gm$ on $\CA=U\cap\prt\Gw$ such that
 \begin{equation}\label{Y3}
\lim_{\ge\to 0}\myint{U\cap\Gs_\ge}{}\gz(\gs(x))u(x)dS(x)=\myint{\CA}{}\gz d\gm\qq\forall \gz\in C_c(U\cap\Gw).
\end{equation}
\es

Notice that any continuous solution of $(\ref{F1})$ in $\Gw$ belongs to  $W^{2,p}_{loc}(\Gw)$ for any $(1\leq p<\infty)$. This previous result yields to a natural definition of the regular boundary points.

\bdef {Trdef}Let $u\in C(\Gw)$ be a positive solution of $(\ref{F1})$. 
A point $z\in\prt\Gw$ is called a regular boundary point for $u$ if there exists an open neighborhood $U$ of $z$ such that $(\ref{F2})$ holds. The set of regular boundary points is a relatively open subset of $\prt\Gw$, denoted by $\CR(u)$. The set $\CS(u)=\prt\Gw\setminus\CR(u)$  is the singular boundary set of $u$. It is a closed set. 
\es

By \rprop{Trlem1} and using a partition of unity, we see that there exists a positive Radon measure $\gm:=\gm_u$
on $\CR(u)$ such that $(\ref{Y3})$ holds with $U$ replaced by $\CR(u)$. The couple $(\gm_u,\CS(u))$ is called the {\bf boundary trace of $u$}.  {\it The main question of the boundary trace problem is to analyse the behaviour of $u$ near the set $\CS(u)$.} \medskip

For any positive good measure $\gm$ on $\prt\Gw$, we denote by $u_\gm$ the solution of $(\ref{D1})$ defined by $(\ref{D9})$-$(\ref{D10})$. 

\bprop{Trlem2} Let  $u\in C(\Gw)\cap W^{2,p}_{loc}(\Gw)$ for any $(1\leq p<\infty)$ be a positive solution of  
$(\ref{F1})$ in $\Gw$ with boundary trace $(\gm_u,\CS(u))$. Then $u\geq u_{\gm_u}$.
\es
\Proof Let $G\subset \prt\Gw$ be a relatively open subset such that $\overline G\subset \CR(u)$
 with a $C^2$ relative boundary $\prt^* G=\overline G\setminus G$. There exists an increasing sequence of 
$C^2$ domains $\Gw_n$ such that $\overline G\subset\prt\Gw_n$, $\prt\Gw_n\setminus \overline G\subset\Gw$ and $\cup_n\Gw_n=\Gw$. For any $n$, let $v:=v_n$ be the solution of 
 \begin{equation}\label{Z1}\left\{\BA {ll}
-\Gd v+Vv=0\qq&\text{in }\Gw_n\\ \phantom{-\Gd v+V}
v=\chi_{_G}\gm\qq&\text{in }\prt\Gw_n.
\EA\right.\end{equation}
Let $u_n$ be the restriction of $u$ to $\Gw_n$. Since $u\in C(\Gw)$ and $Vu\gr\in L^1(\Gw_n)$, there also holds
$Vu\gr_n\in L^1(\Gw_n)$ where we have denoted by $\gr_n$ the first eigenfunction of $-\Gd$ in $W^{1,2}_0(\Gw_n)$. Consequently $u_n$ admits a regular boundary trace $\gm_n$ on $\prt\Gw_n$ (i.e. $\CR(u_n)=\prt\Gw_n$) and $u_n$ is the solution of 
 \begin{equation}\label{Z1}\left\{\BA {ll}
-\Gd v+Vv=0\qq&\text{in }\Gw_n\\ \phantom{-\Gd v+V}
v=\gm_n\qq&\text{in }\prt\Gw_n.
\EA\right.\end{equation}
Furthermore $\gm_n|_G=\chi_{_G}\gm_{u}$. It follows from Brezis estimates and in particular $(\ref{brez3})$ that 
$u_n\leq u$ in $\Gw_n$. Since  $\Gw_n\subset\Gw_{n+1}$, $v_n\leq v_{n+1}$. Moreover 
$$v_n+\BBG^{\Gw_n}[Vv_n]=\BBK^{\Gw_n}[\chi_{_G}\gm]\qquad\text{in }\Gw_n.
$$
Since $\BBK^{\Gw_n}[\chi_{_G}\gm_{u}]\to \BBK^{\Gw}[\chi_{_G}\gm_{u}]$, and the Green kernels $G^{\Gw_n}(x,y)$ are increasing with $n$, it follows from monotone convergence that $v_n\uparrow v$ and there holds
$$v+\BBG^{\Gw}[Vv]=\BBK^{\Gw}[\chi_{_G}\gm_{u}]\qquad\text{in }\Gw.
$$
Thus $v=u_{\chi_{_G}\gm_{u}}$ and $u_{\chi_{_G}\gm_{u}}\leq u$. We can now replace $G$ by a sequence $\{G_k\}$ of relatively open sets with the same properties as $G$, $\overline G_k\subset G_k$ and $\cup_k G_k=\CR(u)$. Then $\{u_{\chi_{_{G_k}}\gm_{u}}\}$ is increasing and converges to some $\tilde u$. Since
$$u_{\chi_{_G{_k}}\gm_{u}}+\BBG^{\Gw}[Vu_{\chi_{_G{_k}}\gm_{u}}]=\BBK^{\Gw}[\chi_{_G{_k}}\gm_{u}],
$$
and $\BBK^{\Gw}[\chi_{_G{_k}}\gm]\uparrow \BBK^{\Gw}[\gm_{u}]$, we derive
$$\tilde u+\BBG^{\Gw}[V\tilde u]=\BBK^{\Gw}[\gm_{u}].
$$
This implies that $\tilde u=u_{\gm_{u} }\leq u$.\qeda

\subsection{The singular part}
The following result is essentially proved in \cite[Lemma 2.8]{MV3}.

\bprop{Trlem3} Let $u\in C(\Gw)$  for any $(1\leq p<\infty)$ be a positive solution of $(\ref{F1})$ and suppose that $z\in \CS(u)$ and that there exists an open neighborhood $U_0$ of $z$ such that 
$u\in L^1(\Gw\cap U_0)$. Then for any open neighborhood $U$ of $z$, there holds
 \begin{equation}\label{Y4}
\lim_{\ge\to 0}\myint{U\cap\Gs_\ge}{}\gz(\gs(x))u(x)dS(x)=\infty.
\end{equation}
\es

As immediate consequences, we have

\bcor {tot1}Assume $u$ satisfies the regularity assumption of \rprop{Trlem2}.  Then for any $z\in\CS(u)$ and any open neighborhood $U$ of $z$, there holds
 \begin{equation}\label{Y5}
\limsup_{\ge\to 0}\myint{U\cap\Gs_\ge}{}\gz(\gs(x))u(x)dS(x)=\infty.
\end{equation}
\es

\bcor {tot2}Assume $u$ satisfies the regularity assumption of \rprop{Trlem2}. If $u\in L^1(\Gw)$,  Then for any $z\in\CS(u)$ and any open neighborhood $U$ of $z$, $(\ref{Y4})$ holds.
\es

The two next results give conditions on $V$ which imply that $\CS(u)=\emptyset$.

\bth{regN=2} Assume $N=2$, $V$ is nonnegative and satisfies $(\ref{cond'})$. If $u$ is a positive solution of  $(\ref{F1})$, then  $\CR(u)=\prt\Gw$.
\es
\Proof We assume that 
\bel{R1}
\myint{\Gw}{}V\gr udx=\infty.
\ee
If $0<\ge\leq \ge_0$, we denote by $(\gr_\ge,\gl_\ge)$ are the normalized first eigenfunction and first eigenvalue of $-\Gd$ in $W^{1,2}_0(\Gw_\ge)$, then
\bel{R2}\lim_{\ge\to 0}\myint{\Gw_\ge}{}V\gr_\ge udx=\infty.
\ee
Because
$$\myint{\Gw_\ge}{}(\gl_\ge+\gr_\ge V)u dx=-\myint{\prt\Gw_\ge}{}\myfrac{\prt\gr_\ge}{\prt{\bf n}}udS,
$$
and
$$c^{-1}\leq-\myfrac{\prt\gr_\ge}{\prt{\bf n}}\leq c,
$$
for some $c>1$ independent of $\ge$, there holds
\bel{R3}\lim_{\ge\to 0}\myint{\prt\Gw_\ge}{}udS=\infty.
\ee
Denote by $m_\ge$ this last integral and set $v_\ge=m^{-1}_\ge u$ and $\gm_\ge=m^{-1}_\ge u|_{\prt\Gw_\ge}$.
Then
\bel{R3'}v_\ge+\BBG^{\Gw_\ge}[Vv_\ge]=\BBK^{\Gw_\ge}[\gm_\ge]\qquad\text{in }\Gw_\ge
\ee
where 
\bel{R4}\BBK^{\Gw_\ge}[\gm_\ge](x)=\myint{\prt\Gw_\ge}{}K^{\Gw_\ge}(x,y) \gm_\ge(y)dS(y)\ee
 is the Poisson potential of $\gm_\ge$ in $\Gw_\ge$ and
$$\BBG^{\Gw_\ge}[Vu](x)=\myint{\Gw_\ge}{}G^{\Gw_\ge}(x,y)V(y)u(y)dy,
$$
the Green potential of $Vu$ in $\Gw_\ge$. Furthermore
\bel{R5}\left\{\BA {ll}
-\Gd v_\ge+Vv_\ge=0\qquad&\text{in }\Gw_\ge\\\phantom{-\Gd v_\ge+V}
v_\ge=\gm_\ge\qquad&\text{in }\prt\Gw_\ge.
\EA\right.\ee
By Brezis estimates and regularity theory for elliptic equations, $\{\chi_{_{\Gw_\ge}}v_\ge\}$ is relatively compact in $L^1(\Gw)$ and in the local uniform topology of $\Gw_{\ge}$. Up to a subsequence $\{\ge_n\}$, $\gm_{\ge_n}$ converges to a probability measure $\gm$ on $\prt\Gw$ in the weak*-topology. It is classical that
$$\BBK^{\Gw_{\ge_n}}[\gm_{\ge_n}]\to \BBK[\gm]
$$
locally uniformly in $\Gw$, and $\chi_{_{\Gw_{\ge_n}}}v_{\ge_n}\to v$  in the local uniform topology of  $\Gw$,  and a.e. in $\Gw$. Because $G^{\Gw_\ge}(x,y)\uparrow G^{\Gw}(x,y)$, there holds for any $x\in \Gw$
\bel{R6}\lim_{n\to\infty}\chi_{_{\Gw_{\ge_n}}}(y)G^{\Gw_{\ge_n}}(x,y)V(y)v_{\ge_n}(y)=G^{\Gw}(x,y)V(y)v(y) \quad\text{for almost all  }y\in \Gw
\ee
Furthermore $v_{\ge_n}\leq \BBK^{\Gw_{\ge_n}}[\gm_{\ge_n}]$ reads
$$v_{\ge_n}(y)\leq c\gr_{\ge_n}(y)
\myint {\prt\Gw_n}{}\myfrac{\gm_{\ge_n}(z)dS(z)}{|y-z|^2}.
$$ 
In order to go to the limit in the expression
\bel{R7}L_n:=\BBG^{\Gw_{\ge_n}}[Vv_{\ge_n}](x)=\myint{\Gw}{}\chi_{_{\Gw_{\ge_n}}}(y)G^{\Gw_{\ge_n}}(x,y)V(y)v_{\ge_n}(y)dy,
\ee
we may assume that $x\in \Gw_{\ge_1}$ where $0<\ge_1\leq\ge_0$ is fixed and write 
$\Gw=\Gw_{\ge_1}\cup\Gw'_{\ge_1}$ where 
$$\Gw'_{\ge_1}=\Gw\setminus\Gw_{\ge_1}:=\{x\in\Gw:\dist (x,\prt\Gw)\leq\ge_1\}$$
and $L_n=M_n+P_n$ where
\bel{R8}M_n=\myint{\Gw_{\ge_1}}{}\chi_{_{\Gw_{\ge_n}}}(y)G^{\Gw_{\ge_n}}(x,y)V(y)v_{\ge_n}(y)dy
\ee
and
\bel{R9}P_n=\myint{\Gw'_{\ge_1}}{}\chi_{_{\Gw_{\ge_n}}}(y)G^{\Gw_{\ge_n}}(x,y)V(y)v_{\ge_n}(y)dy.
\ee
Since
$$\BA {l}
\chi_{_{\Gw_{\ge_1}}}(y)G^{\Gw_{\ge_n}}(x,y)V(y)v_{\ge_n}(y)
\leq c\chi_{_{\Gw_{\ge_1}}}(y)\abs{\ln(|x-y|)}V(y)v_{\ge_n}(y)\\[2mm]\phantom{\chi_{_{\Gw_{\ge_1}}}(y)G^{\Gw_{\ge_n}}(x,y)V(y)v_{\ge_n}(y)}
\leq c\norm{V}_{L^\infty(\Gw_{\ge_1})}\chi_{_{\Gw_{\ge_1}}}(y)\abs{\ln(|x-y|)}v_{\ge_n}(y),
\EA$$
it follows by the dominated convergence theorem that
\bel{R10}\lim_{n\to\infty}M_n=\myint{\Gw_{\ge_1}}{}G^{\Gw}(x,y)V(y)v(y)dy.
\ee
Let $E\subset\Gw$ be a Borel subset. Then $G^{\Gw_{\ge_n}}(x,y)\leq c(x)\gr_{\ge_n}(y)$ if $y\in \Gw'_{\ge_1}$. By Fubini,
\bel{R11}\BA {l}
\myint{\Gw'_{\ge_1}\cap E}{}\chi_{_{\Gw_{\ge_n}}}(y)G^{\Gw_{\ge_n}}(x,y)V(y)v_{\ge_n}(y)dy
\leq cc(x)\myint {\prt\Gw_n}{}\!\!\!\left(\myint{\Gw'_{\ge_1}\cap E}{}\chi_{_{\Gw_{\ge_n}}}\!\!(y)\myfrac{\gr^2_{\ge_n}(y)V(y)}{|y-z|^2}dy\right) \gm_{\ge_n}(z)dS(z)\\[4mm]
\phantom{\myint{A_{\ge_1}\cap E}{}\chi_{_{\Gw_{\ge_n}}}(y)G^{\Gw_{\ge_n}}(x,y)V(y)v_{\ge_n}(y)dy}
\leq cc(x)\displaystyle\max_{z\in\prt\Gw_{\ge_n}}\myint{\Gw'_{\ge_1}\cap E}{}\chi_{_{\Gw_{\ge_n}}}(y)\myfrac{\gr^2_{\ge_n}(y)V(y)}{|y-z|^2}dy
\EA\ee
If $y\in \Gw_{\ge_n}\cap E$, there holds $\gr(y)=\gr_{\ge_n}(y)+\ge_n$. If $z\in  \prt\Gw_{\ge_n}\cap E$ and we denote by $\gs(z)$ the projection of $z$ onto $\prt\Gw$, there holds $|y-\gs(z)|\leq |y-z|+\ge_n$. By monotonicity
\begin{equation}\label{ineq}
\myfrac{\gr_{\ge_n}(y)}{|y-z|}\leq \myfrac{\gr_{\ge_n}(y)+\ge_n}{|y-z|+\ge_n}\leq 
\myfrac{\gr(y)}{|y-\gs(z)|},
\end{equation}
thus
\begin{equation}\label{ineq2}\BA {l}
\myint{\Gw'_{\ge_1}\cap E}{}\chi_{_{\Gw_{\ge_n}}}(y)G^{\Gw_{\ge_n}}(x,y)V(y)v_{\ge_n}(y)dy
\leq cc(x)\displaystyle\max_{z\in\prt\Gw}\myint{\Gw'_{\ge_1}\cap E}{}\chi_{_{\Gw_{\ge_n}}}(y)\myfrac{\gr^2(y)V(y)}{|y-z|^2}dy.
\EA\end{equation}
By $(\ref{cond'})$ this last integral goes to zero if $\abs{\Gw'_{\ge_1}\cap E\cap \Gw_{\ge_n}}\to 0 $. Thus by Vitali's theorem, the sequence of functions $\{\chi_{_{\Gw_{\ge_n}}}(.)G^{\Gw_{\ge_n}}(x,.)V(y)v_{\ge_n}(.)\}_{n\in\BBN}$ is uniformly integrable in $y$, for any $x\in\Gw$. It implies that
\begin{equation}\label{ineq3}
\lim_{n\to\infty}\myint{\Gw}{}\chi_{_{\Gw_{\ge_n}}}(y)G^{\Gw_{\ge_n}}(x,y)V(y)v_{\ge_n}(y)dy
=\myint{\Gw}{}G^\Gw(x,y)V(y)v(y)dy,
\end{equation}
and there holds $v+\BBG[Vv]=\BBK[\gm]$. Since $u=m_\ge v_\ge$ in $\Gw$ and $m_\ge\to\infty$, we get a contradiction since it would imply $u\equiv\infty$. 
\qeda\medskip

In order to deal with the case $N\geq 3$ we introduce an additionnal assumption of stability.

\bth{regN>2} Assume $N\geq 3$. Let $V\in L^\infty_{loc}(\Gw)$, $V\geq 0$ such that
\bel{3-1}
\lim_{\tiny\BA {l}E\text{ Borel}\\|E|\to 0\EA}\myint{E}{}V(y)\myfrac{(\gr(y)-\ge)_+^2}{|y-z|^{N}}dy=0\quad\text {uniformly with respect to }z\in\Gs_\ge \text { and }\ge\in (0,\ge_0].
\ee
 If $u$ is a positive solution of  $(\ref{F1})$, then  $\CR(u)=\prt\Gw$.
\es
\Proof We proceed as in \rth{regN=2}. All the relations $(\ref{R1})$-$(\ref{R10})$ are valid and $(\ref{R11})$ has to be replaced by 
\bel{3-2}\BA {l}
\myint{\Gw'_{\ge_1}\cap E}{}\chi_{_{\Gw_{\ge_n}}}(y)G^{\Gw_{\ge_n}}(x,y)V(y)v_{\ge_n}(y)dy
\leq cc(x)\displaystyle\max_{z\in\Gs_{\ge_n}}\myint{\Gw'_{\ge_1}\cap E}{}\chi_{_{\Gw_{\ge_n}}}(y)\myfrac{\gr^2_{\ge_n}(y)V(y)}{|y-z|^{N+1}}dy.
\EA\ee
Since $(\ref{ineq})$ is no longer valid, $(\ref{ineq})$ is replaced by
\begin{equation}\label{3-3}\BA {l}
\myint{\Gw'_{\ge_1}\cap E}{}\chi_{_{\Gw_{\ge_n}}}(y)G^{\Gw_{\ge_n}}(x,y)V(y)v_{\ge_n}(y)dy
\leq cc(x)\displaystyle\max_{z\in\Gs_{\ge_n}}\myint{ E}{}V(y)\myfrac{(\gr(y)-\ge_n)_+^2}{|y-z|^{N+1}}dy.
\EA\end{equation}
By $(\ref{3-1})$ the left-hand side of $(\ref{3-3})$ goes to zero when $|E|\to 0$, uniformly with respect to $\ge_n$.
This implies that $(\ref{ineq3})$ is still valid and the conclusion of the proof is as in \rth{regN=2}.
\qeda
\medskip

\noindent\Remark A simpler statement which implies $(\ref{3-1})$ is the following.
\begin{equation}\label{3-4}
\lim_{\gd\to 0}\myint{0}{\gd}\left(\myint{B_r(z)}{}V(y)(\gr(y)-\ge)_+^2dy\right)\myfrac{dr}{r^{N+1}}=0,
\end{equation}
uniformly with respect to $0<\ge\leq \ge_0$ and to $z\in\Gs_\ge$. The proof is similar to the one of \rprop{CS}.
\medskip

\noindent\Remark When the function $V$ depends essentially of the distance to $\prt\Gw$ in the sense that
\begin{equation}\label{ineq3}\BA {l}
\abs {V(x)}\leq v(\gr(x))\qquad\forall x\in\Gw,
\EA\end{equation} 
and $v$ satisfies
\begin{equation}\label{ineq4}\BA {l}
\myint{0}{a}tv(t)dt<\infty,
\EA\end{equation} 
Marcus and V\'eron proved \cite[Lemma 7.4]{MV3} that $\CR(u)=\prt\Gw$, for any positive solution  $u$ of $(\ref{F1})$.  This assumption implies also $(\ref{3-1})$.  The proof is similar to the one of \rprop{Tr}.

\subsection{The sweeping method}
 
This method introduced in \cite{RV} for analyzing isolated singulariities of solutions of semilinear equations has been adapted in  \cite{MV0} and \cite{MV4} for defining an extended trace of positive solutions of differential inequalities in particular in the super-critical case. Since the boundary  trace of a positive solutions of $(\ref {F1})$ is known on $\CR(u)$ we shall study the sweeping with measure concentrated on the singular set $\CS(u)$

\bprop{inf} Let $u\in C(\Gw)$ be a positive solution of $(\ref {F1})$ with singular boundary set $\CS(u)$. If  $\gm\in\frak M_+(\CS(u))$ we denote $v_\gm=\inf\{u,u_\gm\}$. Then
 \begin{equation}\label{F2}
-\Gd v_\gm+V(x)v_\gm\geq 0\qquad\text{in }\Gw,
\end{equation}
and $v_\gm$ admits a boundary trace $\gg_u(\gm)\in \frak M_+(\CS(u))$. The mapping $\gm\mapsto \gg_u(\gm)$ is nondecreasing and $\gg_u(\gm)\leq\gm$.
\es
\Proof We know that $(\ref{F2})$ holds
But $Vu_\gm\in L^1_\gr(\Gw)\Longrightarrow Vv_\gm\in L^1_\gr(\Gw)$, if we set $w:=\BBG[Vv_\gm]$, then
$v_\gm+w$ is nonegative and super-harmonic, thus it admits a boundary trace in $\frak M_+(\prt\Gw) $ that we denote by $\gg_u(\gm)$. Clearly $\gg_u(\gm)\leq \gm$ since $v_\gm\leq u_\gm$  and $\gg_u(\gm)$ is nondeacreasing with $\gm$ as $\mu\mapsto u_\gm$ is. Finally, since $v_\gm$ is a supersolution, it is larger that the solution of $(\ref{F1})$ with the same boundary trace $\gg_u(\gm)$, and there holds
 \begin{equation}\label{F3}
u_{\gg_u(\gm)}\leq v_\gm.
\end{equation}

\bprop{Borel} Let
 \begin{equation}\label{F4}
\gn_{_S}(u):=\sup\{\gg_u(\gm):\gm\in   \frak M_+(\CS(u))\}.
\end{equation}
Then $\gn_{_S}(u)$ is a Borel measure on $\CS(u)$. 
\es
\Proof We borrow the proof to Marcus-V\'eron \cite{MV4}, and we naturally extend any positive Radon measure to a positive bounded and regular Borel measure by using the same notation. It is clear that $\gn_{_S}(u):=\gn_{_S}$ is an outer measure in the sense 
 that 
 \begin {eqnarray}\label {borel subadd}\gn_{_S}(\emptyset)=0,\;\mbox { and }\gn_{_S} (A)\leq\sum_{k=1}^\infty\nu (A_{k}),\;\mbox { whenever 
 }A\subset\bigcup_{k=1}^\infty A_{k}.\end {eqnarray}
 Let $A$ and $B\subset \CS(u)$ be disjoint Borel subsets. In 
order to prove that
\begin {eqnarray}\label {borel add}
\gn_{_S} (A\cup B)=\gn_{_S} (A)+\gn_{_S} (B),
\end {eqnarray}
we first notice that the relation holds if $\max \{\gn_{_S} (A),\gn_{_S} 
(B)\}=\infty$. Therefore we assume that $\gn_{_S} (A)$ and $\gn_{_S} (B)$ are 
finite. For $\varepsilon>0$ there exist two bounded positive
measures $\mu_{1}$ and $\mu_{2}$ such that 
$$\gamma_{u}(\mu_{1})(A)\leq \nu (A)\leq \gamma_{u}(\mu_{1})(A)+\varepsilon/2
$$
and
$$\gamma_{u}(\mu_{2})(B)\leq \nu (B)\leq \gamma_{u}(\mu_{2})(B)+\varepsilon/2
$$
Hence
$$\begin {array}{l}
\gn_{_S}(A)+\gn_{_S} (B)\leq \gamma_{u}(\mu_{1})(A)+
\gamma_{u}(\mu_{2})(B)+\varepsilon\\
\phantom {\gn_{_S}(A)+\gn_{_S} (B)}
\leq\gamma_{u}(\mu_{1}+\mu_{2})(A)+\gamma_{u}(\mu_{1}+\mu_{2})(B)
+\varepsilon\\
\phantom {\gn_{_S}(A)+\gn_{_S} (B)}
= \gamma_{u}(\mu_{1}+\mu_{2})(A\cup B)+\varepsilon\\
\phantom {\gn_{_S}(A)+\gn_{_S} (B)}
\leq \gn_{_S} (A\cup B)+\varepsilon.
\end {array}
$$
Therefore $\gn_{_S}$ is a finitely additive measure. If $\{A_{k}\}$ ($k\in 
\BBN$) is a sequence of of disjoint Borel sets and $A=\cup A_{k}$, then
$$\gn_{_S} (A)\geq \gn_{_S} \left(\bigcup_{1\leq k\leq n}A_{k}\right)=\sum_{k=1}^n\gn_{_S} 
(A_{k})\Longrightarrow \gn_{_S} (A)\geq \sum_{k=1}^\infty\gn_{_S} (A_{k}).
$$
By $(\ref {borel subadd})$, it implies that $\gn_{_S}$ is a countably additive measure. \qeda

\bdef{extend} The Borel measure $\gn(u)$ defined by 
 \begin{equation}\label{ext-tr}
\gn(u)(A):=\gn_{_S}(A\cap\CS(u))+\gm_u(A\cap\CR(u)),\qq\forall A\subset\prt\Gw,\, A\text{ Borel},
\end{equation}
is called the extended boundary trace of $u$, denoted by $Tr^{e}(u)$.
\es

\bprop {loc}If $A\subset \CS(u)$ is a Borel set, then
 \begin{equation}\label{rest}
\gn_{_S}(A):=\sup\{\gg_u(\gm)(A):\gm\in  \frak M_+(A)\}.
\end{equation}
\es
\Proof If $\gl,\gl'\in\mathfrak M_+(\CS(u))$  
$$\inf\{u,u_{\gl+\gl'}\}=\inf\{u,u_{\gl}+u_{\gl'}\}\leq \inf\{u,u_{\gl}\}+\inf\{u,u_{\gl'}\}.$$
Since the three above functions admit a boundary trace, it follows that
$$\gg_u(\gl+\gl')\leq \gg_u(\gl)+\gg_u(\gl').
$$
If $A$ is a Borel subset of $\CS(u)$, then $\gm=\gm_{A}+\gm_{A^c}$  where $\gm_{A}=\chi_{_E}\gm$. Thus
$$\gg_u(\gm)\leq \gg_u(\gm_{A})+\gg_u(\gm_{A^c}),
$$
and
$$\gg_u(\gm)(A)\leq \gg_u(\gm_{A})(A)+\gg_u(\gm_{A^c})(A).
$$
Since $\gg_u(\gm_{A^c})\leq \gm_{A^c}$ and $\gm_{A^c}(A)=0$, it follows
$$\gg_u(\gm)(A)\leq \gg_u(\gm_{A})(A).
$$
But $\gm_{A}\leq \gm$, thus $\gg_u(\gm_{A})\leq \gg_u(\gm)$ and finally
 \begin{equation}\label{rest1}
\gg_u(\gm)(A)= \gg_u(\gm_{A})(A).
\end{equation}
If $\gm\in  \frak M_+(A)$, $\gm=\gm_{A}$, thus $(\ref{rest})$ follows.\qeda

\bprop{vanish} There always holds 
 \begin{equation}\label{F5}
\gn(u)(Z^*_V)=0,
\end{equation}
where $Z^*_V$ is the vanishing set of $K_V^\Gw(x,.)$ defined by $(\ref{van1})$.
\es
\Proof This follows from the fact that for any $\gm\in\frak M_+(\prt\Gw)$ concentrated on $Z^*_V$, $u_\gm=0$. Thus $\gg_u(\gm)=0$. If $\gm$ is a general measure, we can write $\gm=\chi_{_{Z^*_V}}\gm+\chi_{_{(Z^*_V)^c}}\gm$, thus $u_\gm=u_{\chi_{_{(Z^*_V)^c}}\gm}$. Because of $(\ref{F3})$ 
$$\gamma_u(\gm)(Z^*_V)=\gamma_u(\chi_{_{(Z^*_V)^c}}\gm)(Z^*_V)\leq (\chi_{_{(Z^*_V)^c}}\gm)(Z^*_V)=0,
$$
thus $(\ref{F5})$ holds.\qeda\medskip

\noindent {\Remark} This process for determining the boundary trace is ineffective if there exist positive solutions $u$ in $\Gw$ such that
$$\lim_{d(x)\to 0}u(x)=\infty.
$$
This is the case if $\Gw=B_R$ and $V(x)=c(R-\abs x)^{-2}$ ($c>0$). In this case $K^\Gw_V(x,.)\equiv 0$. For any $a>0$, there exists a radial solution of
\begin{equation}\label{F6'}
-\Gd u+\myfrac{cu}{(R-|x|)^2}=0\qquad \text {in } B_R
\end{equation}
under the form 
 \begin{equation}\label{F6}
u(r)=u_a(r)=a+c\myint{0}{r}s^{1-N}\myint{0}{s}u(t)\myfrac {t^{N-1}dt}{(R-t)^2}.
\end{equation}
Such a solution is easily obtained by  fixed point, $u(0)=a$ and the above formula shows that $u_a$ blows up when $r\uparrow R$. We do not know if there a exist non-radial positive solutions of $(\ref{F6'})$.
More generaly, if $\Gw$ is a smooth bounded domain, we do not know if there exists a non trivial positive solution of
\begin{equation}\label{F6''}
-\Gd u+\myfrac{c}{d^{2}(x)}u=0\qquad \text {in } \Gw.
\end{equation}

\bth{regN<2} Assume $V\geq 0$ and satisfies $(\ref{cond'})$. If $u$ is a positive solution of $(\ref{F1})$, then  $Tr^{e}(u)=\gn(u)$ is a bounded measure.\es
\Proof Set $\gn=\gn(u)$ and asssume $\gn(\prt\Gw)=\infty$. By dichotomy there exists a decreasing sequence of relatively open domains $D_n\subset\prt\Gw$ such that $\overline D_n\subset D_{n-1}$, diam$\,D_n=r_n\to 0$ as $n\to\infty$, and $\gn(D_n)=\infty$. For each $n$, there exists a Radon measure $\gm_n\in \frak M_+(D_n)$ such that  $\gg_u(\gm_n)(D_n)=n$, 
and 
$$u\geq v_{\gm_n}=\inf\{u,u_{\gm_n}\}\geq u_{\gg_u(\gm_n)}.
$$
Set $m_n=n^{-1}\gg_u(\gm_n)$, then $m_n\in \frak M_+(D_n)$ has total mass $1$ and it converges in the weak*-topology to $\gd_a$, where $\{a\}=\cap_nD_n$. By \rth{stab}, $u_{m_n}$ converges to 
$u_{\gd_a}$. Since $u\geq nu_{m_n}$, it follows that 
$$u\geq \lim_{n\to\infty}nu_{m_n}=\infty,$$
a contradiction. Thus $\gn$ is a bounded Borel measure (and thus outer regular) and it corresponds to a unique Radon measure.\qeda
\medskip

\noindent\Remark If $N=2$, it follows from \rth{regN=2} that $u=u_\gn$ and thus the extended boundary trace coincides with the usual boundary trace. The same property holds if $N\geq 3$, if $(\ref{3-1})$ holds.



\end {document}